\documentclass{article}
\usepackage{amsmath}
\usepackage{amsfonts}
\usepackage{amssymb}

\newtheorem{theorem}{Theorem}

\newtheorem{claim}[theorem]{Claim}
\newtheorem{corollary}[theorem]{Corollary}
\newtheorem{definition}[theorem]{Definition}
\newtheorem{lemma}[theorem]{Lemma}
\newtheorem{proposition}[theorem]{Proposition}
\newtheorem{remark}[theorem]{Remark}
\newenvironment{proof}[1][Proof]{\textbf{#1.} }{\ \rule{0.5em}{0.5em}}

\begin{document}

\title{On a stochastic version of Prouse model in fluid dynamics}
\author{B. Ferrario\\
         Dipartimento di Matematica - Universit\`a di Pavia
        \and F. Flandoli\\
        Dipartimento di Matematica Applicata - Universit\`a di Pisa
  }
\date{ }
\maketitle

\begin{abstract}
A stochastic version of a modified Navier--Stokes equation 
(introduced by Prouse)
is considered in a 3-dimensional torus. 
For equation \eqref{SNS}, we prove existence and uniqueness of
martingale solutions. A different model with the non linearity 
$\Phi(u)=\nu |u|^4u$ is analyzed;
for the structure function of this model, some insights towards an expression
similar to that obtained by the Kolmogorov 1941 theory of turbulence
are presented.
\\
{\bf Key words:} stochastic hydrodynamics, existence and 
uniqueness of martingale solutions, stationary solutions,
structure function in turbulence.
\\
{\bf AMS Subject Classification (2000):} 
 76F55, 
 76M35, 
 76D06, 
 76D03, 
 35Q35. 
\smallskip\\
\end{abstract}

\section{Introduction}

The three dimensional Navier-Stokes equations are a not yet completely
understood mathematical problem, in the sense that there is no proof
of uniqueness of  solutions in the spaces where existence is proved.
This mathematical problem has been investigated since long, also in
connection with
the analysis of how good are the Navier--Stokes equations to model
turbulence.
Some attempts have been made to overcome the problem of uniqueness,
introducing some modification in the Navier-Stokes equations. In this
paper we are concerned with the model proposed by Prouse in
\cite{prouse}.
Here, we study a stochastic version of this problem, as explained below.
As soon as a stochastic equation is introduced, statistical properties
typical of turbulence can be investigated.

We remind that nowadays 
there are many results on stochastic three dimensional
 Navier--Stokes equations
(see, among the others, \cite{bens}, \cite{bcf},\cite{bp},
\cite{cc}, \cite{cf}, \cite{dd}, \cite{f1}, \cite{fg}, \cite{mr}, \cite{vf}); 
however, the uniqueness problem is not solved
also in the stochastic framework.

Let us consider the partial differential equations of
Navier-Stokes type
\begin{equation}\label{SNS}
\left\{
\begin{array}{l}
 du+\left[ - \triangle\Phi(u)+\left(  u\cdot\nabla\right)  u+\nabla p
        -\nabla \text{div } \Phi(u)\right]
 dt  
 =
  G\left(  u\right)  \,dw
\\
 \text{div }u   =0
\\
 u|_{t=0}   =u_{0}
\end{array}
\right.
\end{equation}
where, for $t \ge 0$ and $ x \in \mathcal T \subset \mathbb R^3$, 
$u=u(t,x)$ is the velocity vector field, $p=p(t,x)$  the
pressure field;
$\nu>0$ the viscosity coefficient.
$G$ is an operator acting on the noise and on the velocity;
the vector function $\Phi:\mathbb R^3 \to \mathbb R^3$ 
is defined as follows:
\begin{equation} 
 \Phi(u)=\sigma(|u|)u\quad\text{ with }\quad\left[
\begin{array}[c]{l}
  \sigma\in C^{1}([0,\infty))\\
  \sigma(\xi)\geq\nu>0,\;\sigma^{\prime}(\xi)\geq0\\
  a_1\xi^{b-1} \le \sigma(\xi)\le  a_2\xi^{b-1} 
        \text{ when }\xi > K
\end{array}
\right. \label{Ipo-phi}
\end{equation}
where $a_2\ge a_1>0$ and $b\geq4$.\\
$\Phi$ describes the nonlinear relationship between the stress tensor and
 the deformation velocity tensor, as explained in \cite{prouse}.
When this relationship is linear, then 
$\Phi(u)=\nu u$ and \eqref{SNS}
are the usual Navier--Stokes equations for 
an homogeneous incompressible viscous fluid with random forcing term; 
indeed, the first equation becomes 
$$
  du+\left[ - \nu \triangle u+\left(  u\cdot\nabla\right)  u+\nabla p
        \right]
 dt  
=
  G\left(  u\right)  \,dw
$$

For problem \eqref{SNS}-\eqref{Ipo-phi}, in Section \ref{S3} 
we prove a result on existence and uniqueness of 
martingale solutions (Theorem \ref{teoexist1})
and on existence of stationary martingale solutions (Theorem \ref{teostaz}).

In Section \ref{b=5}, another
model with $\Phi(u)=\nu |u|^4u$ in \eqref{SNS}
is investigated; 
analysis on existence of martingale solutions and stationary
martingale solutions is presented (Theorem \ref{teob5}). 
Moreover, introducing a scaling
transformation suggested by turbulence theory, some insights in  the
behaviour of the function structure  of any order $p$ are shown 
(Claim \ref{K41b5}).

Preliminaries are in Section \ref{prelim}; auxiliary results are in the two 
Appendixes.

\section{Notations and preliminaries}        \label{prelim}

Let the spatial domain be a torus, i.e. the spatial variable $x$ belongs
to $\mathcal T=[0,L]^3$ and 
periodic boundary conditions are assumed. 
$\mathbb{L}^2$ is defined as  the space of vector fields
$u:\mathcal{T}\rightarrow\mathbb{R}^{3}$ with $L^{2}\left(  \mathcal{T}
\right)  $-components. For every $\alpha>0$ and $p>1$, 
$\mathbb{W}^{\alpha,p}$ 
is the space of fields $u\in{\mathbb L}^{p}$ 
with components in the Sobolev space 
$W^{\alpha,p}\left(  \mathcal{T}\right)  $.
For $\alpha<0$, $\mathbb{W}^{\alpha,p}$ is the dual space of 
$\mathbb{W}^{-\alpha,p^\prime}$ 
with $\frac 1p+\frac 1{p^\prime}=1$.
Set $\mathbb{H}^{\alpha}=\mathbb W^{\alpha,2}$.

We introduce the classical spaces for the Navier--Stokes equations
(see, e.g., \cite{temam2}).
$\mathcal{D}^{\infty}$ is defined 
as the space of infinitely differentiable
divergence free periodic fields $u$ on $\mathcal{T}$, with zero mean 
($ \int_{\mathcal{T}} u(x)  dx=0$).
\\ 
Let
$H$ be the closure of $\mathcal{D}^{\infty}$ 
in the $\mathbb{L}^{2}$-topology; it is the space of all fields $u\in
\mathbb{L}^{2}$ such that $\mbox{div}\, u=0$,
$u\cdot n$ on the boundary is periodic, 
$\int_{\mathcal{T}}u\left(  x\right)  dx=0$. We endow $H$
with the inner product
\[
 \left\langle u,v\right\rangle_{H}
 =\frac{1}{L^{3}} \int_{\mathcal{T}} u(x)\cdot v(x) dx
\]
and the associated norm $\left| \cdot \right|_{H}$.

Let $V$ (resp. $D(A)$) be the closure of $\mathcal{D}^{\infty}$ in the
$\mathbb{H}^{1}$-topology (resp. $\mathbb{H}^2 $-topology); 
it is the space of divergence
free, zero mean, periodic elements of $\mathbb{H}^{1}$ 
(resp. of $\mathbb{H}^{2}$). The spaces
$V$ and $D(A)$ are dense and compactly embedded in $H$ (Rellich theorem). Due
to the zero mean condition we also have
\[
 \int_{\mathcal T}\left|  Du\left( x\right)\right|^2 dx\geq\lambda
 \int_{\mathcal T}\left|  u\left( x\right)\right|^2 dx
\]
for every $u\in V$, for some positive constant $\lambda$ (Poincar\'{e}
inequality). 
Here $\left|  Du\left(  x\right)  \right|  ^{2}=\sum_{i,j=1}^{3}\left(
\partial_j u_i(x)\right)  ^{2}$
(and $ \partial_j= \frac{\partial \;}{\partial x_j}$).
So we may endow $V$ with the inner product
\[
 \left\langle u,v\right\rangle_{V}
 = 
 \sum_{i,j=1}^{3} \int_{\mathcal{T}} 
     \partial_j u_i(x)\, \partial_j v_i(x) \, dx
\]
and the associated norm $\left\|  \cdot\right\|  _{V}$.

Let $A:D(A)\subset H\rightarrow H$ be the operator $Au=-\triangle u$
(componentwise).  
There is a complete orthonormal system in $H$  made by the eigenvectors 
$h_{k,j}$ of the
operator $A$ ($Ah_{k,j}=\lambda_{k,j}h_{k,j}$). 
Since the spatial domain is the torus, we know the
expressions of these eigenvectors with their eigenvalues.
Indeed, let $k=(k_1,k_2,k_3)$ with integer components, i.e. $k \in
\mathbb Z^3$. We denote by
$\mathbb Z^3_+$ the half space of $\mathbb Z^3$
defined as $=\{k_1>0\}\cup \{k_1=0,k_2>0\}\cup\{k_1=0,k_2=0,k_3>0\}$.
Then for any $k \in \mathbb Z^3_+$, there exist two unit 
vectors $v_{k,1}$ and $v_{k,2}$, orthogonal to each other and belonging to
the plane orthogonal to $k$.
 Then the (four sequences of) eigenvectors are
$$
 h_{k,1}(x)=\frac {\sqrt 2} {L^{3/2}} v_{k,1} \cos(\tfrac{2\pi}L k \cdot x)\,,
 \qquad
  h_{k,2}(x)=\frac {\sqrt 2}  {L^{3/2}} v_{k,2} \cos(\tfrac{2\pi}Lk \cdot x)
$$
$$
 h_{k,3}(x)=\frac {\sqrt 2}  {L^{3/2}} v_{k,1} \sin(\tfrac{2\pi}Lk \cdot x)\,,
 \qquad
  h_{k,4}(x)=\frac {\sqrt 2}  {L^{3/2}} v_{k,2} \sin(\tfrac{2\pi}Lk \cdot x)
$$
with eigenvalues 
$$
  \lambda_{k,1}=\lambda_{k,2}=\lambda_{k,3}=\lambda_{k,4}
  =
  \frac{(2\pi)^2}{L^2}|k|^2
$$ 
for any $k \in \mathbb Z^3_+$.\\
Hence, $H =\text{span} \{ h_{k,j}: j=1,2,3,4
\mbox{ and } k \in \mathbb Z_+^3\} $
and we set
$H_n=\text{span} \{h_{k,j}: j=1,2,3,4
\mbox{ and } k \in \mathbb Z_+^3, |k|\le n\}$; 
moreover, we denote by
$\pi_n$ the projection operator from $H$ (or any subspace, as $V$ or
$D(A)$) onto $H_n$.
The operators $A$ and $\pi_n$ commute.
\\
We may take the Poincar\'{e} constant $\lambda$ above
equal to $(2\pi)^2/L^2$ (the first eigenvalue of $A$). Notice that we have
\[
 \left\langle Au,u\right\rangle _{H}=\left\|  u\right\|  _{V}^{2}
\]
for every $u\in D(A)$, so in particular
\[
 \left\langle Au,u\right\rangle _{H}\geq 
  \frac{(2\pi)^2}{L^2}\left|  u\right|  _{H}^{2}.
\]

Let $V^{\prime}$ be the dual of $V$ with respect to the $H$-norm; 
with proper identifications we have
$V\subset H\subset V^{\prime}$ with continuous injections, and the scalar
product $\left\langle \cdot,\cdot\right\rangle _{H}$ extends to the dual
pairing $\left\langle \cdot,\cdot\right\rangle _{V,V^{\prime}}$ between $V $
and $V^{\prime}$ and to the dual
pairing $\left\langle \cdot,\cdot\right\rangle 
_{\mathbb L^q,\mathbb L^{q^\prime}}$ between 
$\mathbb L^q$ and $\mathbb L^{q^\prime}$ ($\frac 1q
+\frac 1{q^\prime}=1$).

Let $B\left(  \cdot,\cdot\right)  :V\times V\rightarrow V^{\prime}$ be the
bilinear operator defined as
\begin{equation}\label{defB}
 \left\langle w,B\left(  u,v\right)  \right\rangle _{V,V^{\prime}}
=
 \sum_{i,j=1}^{3}\int_{\mathcal T} u_i (\partial_i v_j) w_j dx
\end{equation}
for every $u,v,w\in V$. 
By the incompressibility condition, we have
\begin{equation}\label{incomp}
 \langle B\left(  u,v\right),v  \rangle =0, \qquad
 \langle B\left(  u,v\right),w  \rangle = 
       -\langle B\left(  u,w\right),v  \rangle 
\end{equation}
Using the latter relationship, by H\"older inequality we estimate 
\begin{equation}\label{BconL4}
 |B(u,u)|_{V^\prime}
 =\sup_{\|\psi\|_V\le 1} |\langle B(u,u),\psi\rangle|
 \le |u|^2_{\mathbb L^4}
\end{equation}

We list here a number of inequalities.

\begin{lemma}
                      \label{2}
\[
  \langle A\Phi(u),u\rangle_H \geq \nu\Vert u\Vert_V^{2}
\]
\end{lemma}
\begin{proof} 
We have
\[
  \langle A\Phi(u),u\rangle_H=
  \langle \Phi(u), u\rangle_V=
  \int_{\mathcal T} \sum_{i,k=1}^3\left[\partial_k \Phi_i(u)\right]
     \partial_k u_i \,dx
\]
The estimate on $\sum_{i,k=1}^3\left[\partial_k \Phi_i(u)\right] \partial_k u_i$
comes from \cite{prouse}.
\end{proof}

\begin{lemma}
              \label{3}
\[
 \langle\Phi(u^{(1)})-\Phi(u^{(2)}),u^{(1)}-u^{(2)}\rangle_H
\geq
 \nu |u^{(1)}-u^{(2)}|_H^2
\]
\end{lemma}
\begin{proof}
The proof is by \cite{prouse}. We rewrite it here, because we shall
need it in Section \ref{b=5}.

Set $\sigma(|u|)=\nu+\tilde \sigma(|u|)$ with $\tilde \sigma^\prime\ge
0$.
Then
\begin{align*}
&[\tilde\sigma (|u^{(1)}|) u^{(1)} - \tilde\sigma (|u^{(2)}|) u^{(2)}]\cdot [u^{(1)}-u^{(2)}]
\\
&=
 \tilde\sigma (|u^{(1)}|) |u^{(1)}|^2 + \tilde\sigma (|u^{(2)}|) |u^{(2)}|^2
 - \tilde\sigma (|u^{(1)}|) u^{(1)} \cdot u^{(2)}
 -  \tilde\sigma (|u^{(2)}|) u^{(1)} \cdot u^{(2)}
\\
& 
\ge
  \tilde\sigma (|u^{(1)}|) |u^{(1)}|^2 +\tilde \sigma (|u^{(2)}|) |u^{(2)}|^2
\\
&\hspace*{5mm}
- \frac 12  \tilde\sigma (|u^{(1)}|) |u^{(1)}|^2
- \frac 12  \tilde\sigma (|u^{(2)}|) |u^{(2)}|^2
- \frac 12  \tilde\sigma (|u^{(1)}|) |u^{(2)}|^2
- \frac 12  \tilde\sigma (|u^{(2)}|) |u^{(1)}|^2
\\
&
=
 \frac 12  \tilde\sigma (|u^{(1)}|) |u^{(1)}|^2 
 + \frac 12 \tilde\sigma (|u^{(2)}|) |u^{(2)}|^2
 - \frac 12 \tilde\sigma (|u^{(1)}|) |u^{(2)}|^2
 - \frac 12 \tilde\sigma (|u^{(2)}|) |u^{(1)}|^2
\\
&
=\frac 12 \left[\tilde \sigma (|u^{(1)}|)-\tilde \sigma (|u^{(2)}|)\right]\,
 \left[|u^{(1)}|-|u^{(2)}|\right] \, \left[|u^{(1)}|+|u^{(2)}|\right]
\\
&
\ge 0
\end{align*}
Hence
\begin{align*}
&[\sigma (|u^{(1)}|) u^{(1)} - \sigma (|u^{(2)}|) u^{(2)}]\cdot [u^{(1)}-u^{(2)}]
\\
&=
 \nu |u^{(1)}-u^{(2)}|^2+ 
 [\tilde \sigma (|u^{(1)}|) u^{(1)} - \tilde \sigma (|u^{(2)}|) u^{(2)}]\cdot [u^{(1)}-u^{(2)}]
\\
&
\ge 
 \nu |u^{(1)}-u^{(2)}|^2
\end{align*}
\end{proof}

Next lemma is crucial to prove uniqueness. Notice that the regularity
$u \in L^5(0,T;\mathbb L^5)$
is needed here. The weak solutions of the  Navier--Stokes
equations (deterministic or stochastic), which are known to exist, 
are not proved to have  such a regularity; 
here the modified term $\Phi$
(with $b \ge 4$) plays its role. We remind that Prodi \cite{prodi} proved
uniqueness for the deterministic three dimensional Navier--Stokes equations,
if $u \in L^{\frac {2q}{q-3}}(0,T;\mathbb L^q)$
for some $3<q\le \infty$. For $q=5$ the required regularity is $u \in 
L^5(0,T;\mathbb L^5)$ and this implies uniqueness also in the Prouse
model (see \cite{prouse}).

\begin{lemma} 
                     \label{1} 
If $u\in \mathbb L^5$ and $v\in H$, then for any $\nu>0$
\[
 |\langle B(u,v),A^{-1}\pi_m v\rangle|
 \le
 \frac{\nu}{4}|v|_H^2 + C_B|u|_{\mathbb L^5}^5 |\pi_m v|_{V^{\prime}}^2
\]
\[
 |\langle B(v,u),A^{-1}\pi_m v\rangle|
 \le
 \frac{\nu}{4}|v|_H^2 + C_B|u|_{\mathbb L^5}^5 |\pi_m v|_{V^{\prime}}^2
\]
for some positive constant $C_B$.
\end{lemma}
\begin{proof}
In \cite{prouse}, there is a very similar lemma, but with $v$ instead
of $\pi_m v$ (here we consider any finite projection operator
$\pi_m$). 
Following the lines of that proof, we get our result. 
\end{proof}

         \vspace{2mm}
{\bf Properties of $G$}\\
Let $G:H\rightarrow L\left(  H\right)  $ be a mapping with the properties
\begin{equation}\label{crescG}
 \|G(u)\|^2_{HS(H)}
  \le
 \lambda_0|u|^2_H+\rho
\end{equation}
and
\begin{equation} \label{Glip}
  \left\|A^{-1/2}[G(v)-G(z)]\right\|^2_{HS(H)}
 \le
  L_G |v-z|^2_{V^\prime}
\end{equation}
Here $\left\|  T\right\|  _{HS\left(  H\right)  }$ is the Hilbert-Schmidt
norm of an operator in $H$, defined as
\[
 \left\| T \right\|_{HS(H)}^2 = \sum_{j=1}^4
      \sum _{k \in \mathbb Z_+^3}\left|T h_{k,j}\right|_H^2
\]

Now, we project equation \eqref{SNS} onto the space of 
divergence free vectors fields; 
both the $\nabla$-terms desappear, as 
when we deal with the Navier--Stokes equations (see, e.g.,
\cite{temam}).
Then, we obtain an
evolution equation  (still formally), which with our notations is 
\begin{equation}
 du+[A\Phi(u)+B(u,u)]\,dt=G\left(  u\right)  \,dw,
     \quad u(0)=u_{0}
 \label{abstr SNS}
\end{equation}
From now on, $\Phi$ will be assumed to satisfy \eqref{Ipo-phi} for 
a given $b\ge 4$.
\\
The rigorous interpretation  of this equation
will be given in the sequel, but for the time
being let us at least write it in weak form
\begin{multline} \label{weakSNSeq}
 \left\langle u_t,\psi\right\rangle_H
   +\int_0^t \left\langle\Phi(u_s),A\psi\right\rangle
  _{\mathbb L^{1+\frac 1b},\mathbb L^{1+b}} ds
   -\int_0^t \left\langle B\left(u_s,\psi\right),u_s\right\rangle
     _{\mathbb L^{\frac 43},\mathbb L^4} ds
\\
  =\left\langle u_0,\psi\right\rangle_H
   +\int_0^t\left\langle G(u_s) \,dw_s,\psi\right\rangle_H  
\end{multline}
with $\psi\in\mathcal{D}^{\infty}$ and $0<t<\infty$.
\\
We assume that $w$ is a cylindrical Wiener process in $H$ (see,
e.g., \cite{dpz}). We can represent it as follows.
Suppose we are given  a Brownian stochastic basis, i.e.
a probability space $\left(\mathcal W,\mathcal{F},Q\right)  $,
a filtration $( \mathcal{F}_{t})_{t\geq0}$ and a 
sequence $\{\beta_{k,j}(t)\}_{k,j}$ of independent
Brownian motions on $\left(\mathcal W,\mathcal{F},\left(\mathcal{F}_{t}\right)
_{t\geq0},Q\right)  $. Namely, for $k \in \mathbb Z^3_{+}$ and
$j=1,2,3,4$, 
the real valued processes $\beta_{k,j}(t)$
are independent, adapted to $\left(  \mathcal{F}_{t}\right)  _{t\geq0}$, 
continuous for $t \ge 0$ and null at $t=0$, 
with increments $\beta_{k,j}(t)-\beta_{k,j}(s)$
that are $N\left(  0,t-s \right)  $-distributed
and independent of $\mathcal{F}_{s}$.
Then  
\begin{equation} \label{cilindrico}
 w(t)=\sum_{j=1}^4 \sum_{k\in \mathbb Z^3_+} \beta_{k,j}(t) h_{k,j}
\end{equation} 
is a cylindrical Wiener process in $H$.

The convergence of this series requires
proper distributional topologies. The stochastic integral in equation 
\eqref{weakSNSeq}  is
well defined under the Hilbert-Schmidt assumption made on $G$ (see 
\cite{dpz} for details).

\section{Well posedness} \label{S3}

\subsection{Concepts of solution}

Consider the abstract (formal) stochastic evolution equation (\ref{abstr SNS})
and its weak formulation over test functions (\ref{weakSNSeq}). 
We have
\begin{align*}
 \int_0^t \left|\left\langle\Phi(u_s),A\psi\right\rangle
       _{\mathbb L^{1+\frac 1b}, \mathbb L^{1+b}}\right| ds
&
\le 
   \int_0^t \left|\Phi(u_s)\right|_{\mathbb L^{1+\frac 1b}}
            \left| A \psi \right|_{\mathbb L^{1+b}} ds
\\
&
\le 
 C_\psi \int_0^t (1+\left|u_s\right|^{b}_{\mathbb L^{1+b}})ds                
\end{align*}
because
\begin{equation}\label{1+1b}
\begin{split}
  \left| \Phi(u)  \right|_{\mathbb L^{1+\frac 1b}}^{1+\frac 1b}
 &=
  \int_\mathcal T |\Phi(u(x))|^{1+\frac 1b} 1_{\{|u(x)|\leq K\}}dx+
  \int_\mathcal T |\Phi(u(x))|^{1+\frac 1b} 1_{\{|u(x)|> K\}}dx
\\
 & \stackrel{\text{by }\eqref{Ipo-phi}}{\le}
  K^{1+\frac 1b} \int_\mathcal T |\sigma(|u(x)|)|^{1+\frac 1b}1_{\{|u(x)|\leq K\}} dx+ 
  \int_\mathcal T (a_2 |u(x)|^b)^{1+\frac 1b} dx
\\
 & \leq 
  C_{\Phi} (1+\int_\mathcal T |u(x)|^{1+b} dx )
\end{split}
\end{equation}
since $\sigma \in C^1$ implies that $\sigma$ is bounded on $[0,K]$.
Then, in equation \eqref{weakSNSeq} 
the term $ \int_0^t \left\langle\Phi(u_s),A\psi\right\rangle ds$
is well defined for functions
$u$ that live in $L^{1+b}(0,T;\mathbb L^{1+b})$, $T>0$. 
\\
Moreover, 
\begin{equation}\label{senso}
   \int_0^t \left|\langle B(u_s,\psi),u_s\rangle
        _{\mathbb L^{\frac 43},\mathbb L^4}  \right|  ds 
 \leq 
   \int_0^t \left|  u_s\right|_{\mathbb L^4}^2 \left\|\psi\right\|_V ds
 \leq 
      C_{\psi} \int_0^t \left|  u_s\right|_{\mathbb L^4}^2 ds
\end{equation}
Hence, in equation \eqref{weakSNSeq}
the  term $\int_0^t \langle B(u_s,\psi),u_s\rangle ds$ 
is well defined for functions
$u$ that live in $L^2(0,T;\mathbb L^4)$.

We conclude, in both cases, that given $b\ge 4$ the regularity 
 $u \in  L^{1+b}(0,T;\mathbb L^{1+b})$ is enough to define these 
quantities. Moreover, 
from now on the duality pairing for these two terms
 has to be understood in the sense above specified (as written also
 in equation \eqref{weakSNSeq}).

As in the deterministic case, strong continuity of trajectories in $H$ is an
open problem. There will be strong continuity in weaker spaces (like
$\mathbb W^{-2-\theta,1+\frac 1b}$), 
and a uniform bound in $H$. Let $H_{\sigma}$ be the space $H$
with the weak topology. Since
\[
 C([0,T];\mathbb W^{-2-\theta,1+\frac 1b})  
  \cap L^{\infty}(0,T;H)  
   \subset C([0,T];H_{\sigma})
\]
then the trajectories of the solutions will be at least weakly continuous
in $H$
(see \cite{temam} pg. 263).

Given a separable Banach space $\mathbb W^\prime$ (it will be
$\mathbb W^\prime=\mathbb W^{-2-\theta,1+\frac 1b}$), let us set
$$
 \Omega=C(  [0,\infty);\mathbb W^\prime) 
$$
and denote by $\left(  \xi_{t}\right)  _{t\geq0}$ the canonical process
($\xi_{t}\left(  \omega\right)  =\omega_{t}$), by $F$ the Borel $\sigma
$-algebra in $\Omega$ and by $F_{t}$ the $\sigma$-algebra generated by the
events $\left(  \xi_{s}\in A\right)  $ with $s\in\left[  0,t\right]  $ and
$A\in\mathcal{B}\left(\mathbb W^\prime\right)  $.

\begin{definition}
[solution to the martingale problem]
Given a probability measure $\mu_{0}$ on
$H$, we say that a probability measure $P$ on $\left(  \Omega,F\right)  $ is a
solution of the martingale problem associated to equation (\ref{abstr SNS}) with
initial law $\mu_{0}$ if

\begin{enumerate}
\item[\lbrack MP1\rbrack] 
for every $T>0$
\[
 P\left(  \sup_{t\in [0,T]} \left|\xi_t\right|_H
  +\int_0^T \left\|\xi_s\right\|_V^2 ds
  +\int_0^T \left|\xi_s\right|_{\mathbb L^{1+b}}^{1+b} ds<\infty\right)  =1
\]

\item[\lbrack MP2\rbrack] 
for every $\psi\in\mathcal{D}^{\infty}$ the process
$M_{t}^{\psi}$ defined $P$-a.s on $\left(  \Omega,F\right)  $ as
$$
 M_t^\psi 
 :=
  \langle \xi_t,\psi\rangle_H
  -\langle \xi_0,\psi\rangle_H
  -\int_0^t \langle\Phi\left(\xi_s\right) ,A\psi\rangle ds
  +\int_0^t \langle B\left(\xi_s,\psi\right),\xi_s \rangle ds
$$
is square integrable and $\left(  M_{t}^{\psi},F_{t},P\right)  $ is a
continuous martingale with quadratic variation
\[
 \left[  M^{\psi}\right]_t
 =
 \int_0^t |G(\xi_s)\psi|_H^2 ds
\]

\item[\lbrack MP3\rbrack] $\mu_{0}=\Pi_{0}P$, where $\Pi_0$ denotes the
  restriction on $F_0$.
\end{enumerate}
\end{definition}

\begin{remark}
A solution of the martingale problem is also a  weak
solution. The definition of  weak solution is as follows:
there exists
a Brownian stochastic basis $\left( \mathcal W,\mathcal{F},\left(  \mathcal{F}
_{t}\right)  _{t\geq0},Q,\left(  \beta_{i}\left(  t\right)  \right)
_{t\geq0;k,j}\right)  $ 
and a $\mathbb W^\prime$-valued process $u$ on
$\left(\mathcal  W,\mathcal{F},Q\right)  $ such that

\begin{enumerate}
\item [\lbrack WM1\rbrack] $u$ is a continuous adapted process 
in $\mathbb W^\prime$ and
\[
u\left(  .,\omega\right)  \in L^{\infty}\left(  0,T;H\right)  \cap
L^{2}\left(0,T;V\right) \cap L^{1+b}(0,T;\mathbb L^{1+b})\quad Q\text{-a.s.}
\]
for every $T>0$

\item[\lbrack WM2\rbrack] (\ref{weakSNSeq}) is satisfied $Q$-a.s.

\item[\lbrack WM3\rbrack] $u(0)$ has law $\mu_{0}$.
\end{enumerate}
\end{remark}

Finally, in this context, we call strong solution 
a process $u$ satisfying
the three above properties 
 on any a priori given stochastic basis.

\subsection{Main result}

\begin{theorem}
\label{teoexist1}Let $\mu$ be a measure on $H$ such that 
$m_p:=\int_H |v|_H^p \mu(dv)  <\infty$ for some $p>2$. 
Then there
exists one and only one solution to the martingale problem 
\eqref{weakSNSeq} with initial
condition $\mu$.

Moreover, two strong 
solutions on the same Brownian stochastic basis coincide a.s. 
\end{theorem}

\begin{proof}
\textbf{Step 1} (Galerkin approximations). Let
\[
\left(\mathcal W, \mathcal{F}, (\mathcal{F}_{t})_{t\geq0},Q,
 \left(\beta_{k,j}(t) \right)_{t\geq0;k,j}\right)
\]
be a Brownian stochastic basis supporting also an $\mathcal{F}_{0}$-measurable
r.v. $u_{0}:W\rightarrow H$ with law $\mu$. 
For every $n$, let $u_{0}^{n}:=\pi_{n}u_{0}$
and consider the Galerkin system
\begin{equation}\label{Galerkin SNS}
 du_t^{n}+[A\Phi(u_t^{n})+\pi_{n}B(u_t^{n},u_t^{n})]\,dt
=
 \pi_{n}G\left(u_t^{n}\right)  \,dw_t,
     \quad u^{n}(0)=u^n_0
\end{equation}
obtained by applying the projection operator $\pi_n$ to both sides of
equation \eqref{abstr SNS}.
(Notice that $\pi_n A \Phi(u^n_t)=A\Phi(u^n_t)$.)
\\
Equation \eqref{Galerkin SNS} is a stochastic
ordinary equation in the finite-dimensional Hilbert space $H_n$.
\\
Local existence and uniqueness (on a random
time interval) 
is classical, since the nonlinearities are locally Lipschitz
continuous (see, e.g., \cite{protter}). 
Global existence  is then a consequence of the a priori estimates
given in Appendix 1.
There, defined $\tau_R^n=\inf\left\{  t\ge 0: |u^n_t|_H^2=R\right\}$
we shall prove that, for any $T>0$
\begin{equation}\label{mediaP}
 E \sup_{0\le t\le T} |u^n_{t\wedge \tau^n_R}|_H^p \le C_1
\end{equation}
\begin{equation} \label{stima-2}
 E \int_0^T \|u^n_{s\wedge\tau_R}\|_V^2 1_{\{s< \tau_R\}}ds
\le C_2
\end{equation}
\begin{equation} \label{stima-3}
 E \int_0^T |u^n_{s\wedge \tau^n_R}|^{1+b}_{\mathbb L^{1+b}}
 1_{\{s<\tau^n_R\}} ds
 \le C_3
\end{equation}
for some positive constants $C_1=C_1(p,T,\lambda_0,\rho,m_p)$,
$C_2=C_2(T,\lambda_0,\rho,m_2)$, $C_3=C_3(a_1,C_1,C_2)$, 
independent of $n$ and $R$.

Now, assume first that the initial velocity is bounded:
$|u_0|_H\le K$.
Take $R>K$; so $\tau_R^n>0$ $Q$-a.s..
The solution $u^n_t$ to the Galerkin system
\eqref{Galerkin SNS} is defined at least in the time interval $[0,\tau_R^n)$.
Since we know from \eqref{mediaP}
that
$$
 E \sup_{t \in [0,T]}|u^n_{t\wedge \tau_R^n}|^2_H \le \tilde C_1
$$
for some constant $\tilde C_1=C_1^{2/p}$ independent of $n$ and $R$,
we have 
$$
 E\left(1_{\{\tau^n_R<T\}} |u^n_{T\wedge \tau_R^n}|_H^2\right)
 \le \tilde C_1
$$
for $T>0$ fixed.
Moreover
$$
 Q(\tau^n_R<T) = E 1_{\{\tau^n_R<T\}}
 =
 \frac 1R E \left(1_{\{\tau^n_R<T\}} |u^n_{T\wedge \tau_R^n}|_H^2\right)
$$
because $|u^n_{T\wedge \tau_R^n}|^2_H=R$ on the set $\{\tau^n_R<T\}$.
Hence
$$
 Q(\tau^n_R<T) 
 \le
 \frac{\tilde C_1}{R}
$$
Notice that 
$\tau^n_{\tilde R}> \tau^n_R$  for $\tilde R>R$.
Therefore, setting $\tau^n_\infty={\displaystyle \sup_{R>K}} \tau^n_R$ the 
process $u^n_t$ is defined  for $t \in [0,\tau^n_\infty)$.
But we have
$$
 Q(\tau^n_\infty<T)\le
 Q(\tau^n_R<T)
 \le
  \frac{\tilde C_1}{R}
                \qquad \forall R
$$
Hence 
$$
  Q(\tau^n_\infty<T)=0
$$
and finally we conclude that $u^n_t$ is a solution for $t \in
[0,T)$. 
Since $T$ has been chosen arbitrarily, we conclude that the Galerkin
solution is defined on any finite time interval.
 
For a general initial velocity satisfying the assumption of Theorem 
\ref{teoexist1}, we proceed as follows. Let ${\mathcal W}_K
\in \mathcal F$ be defined as 
${\mathcal W}_K=\{|u_0|^2_H\le K\}$; we have
$Q(\cup_K {\mathcal W}_K)=Q(|u_0|_H<\infty)=1$. 
Define ${u_0}_K$ as $u_0$ on ${\mathcal W}_K$
and 0 otherwise.
Let ${u^n_t}_K$ be the unique solution to the Galerkin system ($\forall
t \ge 0$) with
initial condition
${u_0}_K$. If $\tilde K> K$, then 
$$
 Q\left\{ {\mathcal W}_K\cap \{{u^n_t}_{\tilde K}={u^n_t}_K \; 
       \forall t\ge 0\}\right\}
 = Q\{ {\mathcal W}_K\}
$$
We may uniquely define a process ${u^n_t}_\infty$ on
$\mathcal W^\prime=\cup_K {\mathcal W}_K$ as ${u^n_t}_\infty={u^n_t}_K$ on
${\mathcal W}_K$.
Looking at the Galerkin equation in the integral form, it is clear
that ${u^n_t}_\infty$ solves the equation on $\mathcal W ^\prime$.
But $Q(\mathcal W^\prime)=1$. Thus we have proved the existence of a
global solution to the Galerkin system for any 
initial velocity with $m_p<\infty$ for some $p>2$.
This solution is a continuous adapted Markov process in $H_n$
(uniqueness holds for the Galerkin problem; it can be checked direclty
or obtained as a byproduct of next Step 5).

Hence we have proved that, for any $T<\infty$
\begin{equation}  \tag{\ref{mediaP}$'$}
  E \sup_{0\le t\le T} |u^n_t|_H^p \le C_1
\end{equation}
\begin{equation}  \tag{\ref{stima-2}$'$}
 E \int_0^T \|u^n_s\|_V^2 ds
\le C_2
\end{equation}
\begin{equation}  \tag{\ref{stima-3}$'$}
 E \int_0^T |u^n_s|^{1+b}_{\mathbb L^{1+b}}
 ds
 \le C_3
\end{equation}

From these estimates, we also get the following one.
Given $\psi \in \mathcal D^\infty$ and $\varepsilon \in (0,2)$, we have 
$$
\begin{array}{rl}
 \displaystyle
 E\left| \int_0^t\langle \pi_nG(u^n_s)dw_s,\psi\rangle\right|^{2+\varepsilon}
 &\le \displaystyle
 \left( E \left| \int_0^t\langle \pi_n G(u^n_s)dw_s,\psi\rangle\right|^4
         \right)^{\frac{2+\varepsilon}4}
\\
  & =  \displaystyle
 \left( C\left( E  \int_0^t |\pi_n G(u^n_s)\psi|^2_Hds\right)^2 
      \right)^{\frac{2+\varepsilon}4}
           \text{ by Gaussianity  }
\\
  & \le  \displaystyle
  C |\psi|_H^{2+\varepsilon}
   \left(  \left( E \int_0^t \|G(u^n_s)\|_{HS(H)}^2 ds\right)^2
   \right)^{\frac{2+\varepsilon}4}
\\
  &\le  \displaystyle
 C |\psi|^{2+\varepsilon}_H
  E \int_0^t(|u^n_s|_H^{2+\varepsilon}+1)ds
\end{array}
$$
Here (and in the following) 
$C$ denotes different positive  constants, independent of $n$.
Taking $2+\varepsilon \le p$ and bearing in mind \eqref{mediaP}, 
we conclude that for any finite $t$
\begin{equation}\label{2e-marting}
 \sup_n E \left| \int_0^t\langle \pi_nG(u^n_s)dw_s,\psi\rangle\right|^{2+\varepsilon}<\infty
\end{equation}
Here the limitation $\varepsilon <2$ can be easily removed, but in the
sequel it will be enough to consider a positive quantity
$\varepsilon $ as small as we want.

\textbf{Step 2} (time regularity and reformulation in path space). 
In view of the time regularity, equation (\ref{Galerkin SNS}) has the
form
\[
u_{t}^{n}=u_{0}^{n}+I_{t}^{n}+J_{t}^{n}+K_{t}^{n}
\]
where
\begin{align*}
& I_{t}^{n}=-\int_{0}^{t} A\Phi\left(  u_{s}^{n}\right) ds
\\
&
 J_{t}^{n}=-\int_{0}^{t} \pi_{n} B\left(  u_{s}^{n},u_{s}^{n}\right)  ds
\\
&
 K_{t}^{n}=\int_{0}^{t}\pi_{n}G\left(  u_{s}^{n}\right)  dw_{s}
\end{align*}
For the first term we have 
\begin{align*}
   \left\|  I_{\cdot}^{n}\right\|_{W^{1,1 + \frac 1b}\left(
    0,T;\mathbb W^{-2,1+\frac 1b}\right)  }^{1+\frac 1b}
 & \leq 
   C \int_0^T \left|  A\Phi\left(  u_{s}^{n}\right)  \right|
   _{\mathbb W^{-2,1+\frac 1b}}^{1+\frac 1b}ds
\\
 & \le
  C \int_0^T \left| \Phi\left(  u_{s}^{n}\right)  \right|
  _{\mathbb L^{1+\frac 1b}}^{1+\frac 1b}ds
\\
 & \leq
    C (T+\int_0^T \left| u_s^n \right|
   _{\mathbb L^{1+b}}^{1+b} ds)
\end{align*}
according to \eqref{1+1b}.
\\
For $J_{t}^{n}$, using \eqref{BconL4} 
we have
\begin{align*}
 \left\|  J_{\cdot}^{n}\right\|_{W^{1,2}\left(0,T;V^\prime \right)}^2  
 & \leq 
  C \int_{0}^{T}\left|  B\left(  u_{s}^{n},u_{s}
        ^{n}\right)  \right|  _{V^\prime}^{2}ds
\\
 & \leq 
  C \int_{0}^{T}\left|  u_{s}^{n}\right|  _{\mathbb L^4}^{2}ds
\\
 & \le
 C (T+\int_0^T \left| u_s^n \right|
   _{\mathbb L^{1+b}}^{1+b} ds)
\end{align*}
Finally, for every $q>1$, $\alpha\in\left(  0,\frac 12\right)  $,
$T>0$, we have (see, e.g., \cite{fg})
\[
 E\left\|  K_{\cdot}^{n}\right\|  _{W^{\alpha,q}\left(  0,T;H\right)
 }^{q}
 \le
 C E\int_0^T \|\pi_n G(u^n_s)\|^q_{HS(H)}ds
\]
and by \eqref{crescG} and the mean estimates of the
previous step
we conclude that
\[
 E\left\|  K_{\cdot}^{n}\right\|  _{W^{\alpha,p}\left(  0,T;H\right)
 }^{p}
\le
 \tilde C
\]
($\tilde  C$ independent of $n$ and $p>2$ 
as stated in Theorem \ref{teoexist1}.)

Therefore,  for $\alpha\in (0,\frac 12)$
\begin{equation}\label{3+reg}
 u^n \in W^{1,1 + \frac 1b}\left(0,T;\mathbb W^{-2,1+\frac 1b}\right)
 +W^{1,2}\left(  0,T;V^\prime \right) 
 +W^{\alpha,p}\left(  0,T;H\right) 
\end{equation}
in mean.
Notice that
$H\subset V^\prime \subset \mathbb W^{-2,1+\frac 1b}$,
$W^{1,2}(0,T)\subset  W^{1,1 + \frac 1b}(0,T) \subset 
 W^{\alpha,1 + \frac 1b}(0,T)$
and
$W^{\alpha,p}(0,T)\subset W^{\alpha,1 + \frac 1b}(0,T)$.
\\
We conclude that, in mean 
$$
 u^n \in  W^{\alpha,1+\frac 1b}(0,T;\mathbb W^{-2,1+\frac 1b})
$$

Under the
embedding $H_{n}\subset H$, we have that $\left(  u_{t}^{n}\right)  _{t\geq0}$
is a continuous adapted process in $H$, so it defines a measure $P_{n}$ on
$C\left(  [0,\infty);H\right)  $, and thus on $\left(
  \Omega,F\right)$. 
Actually, $P_n$ is concentrated on $C([0,\infty);H_n)$.
For every $\alpha\in\left(  0,\frac 12\right)  $,  $T>0$, 
the above estimates may be rewritten as
\[
 E^{P_n} \left[  \sup_{t\in\left[  0,T\right]  }\left|  \xi_{t}\right|^p_{H}
 +\int_{0}^{T}\left\|  \xi_{s}\right\|  _{V}^{2}ds
 + \int_0^T \left| \xi_s \right|_{\mathbb L^{1+b}}^{1+b}ds
 \right]  
\leq
 C_{4}\left(T,\lambda_0,\rho,m_p,b\right)
\]
and
\begin{equation}\label{reg-spaz}
E^{P_n}\left[  \left| \xi\right|_{W^{\alpha,1+\frac
      1b}(0,T;\mathbb W^{-2,1+\frac 1b})}\right]  
  \le 
  C_5 \left(T,\lambda_0,\rho,m_p,b\right) 
\end{equation}
for any $n$.

Relationships $(\ref{mediaP}')$-$(\ref{stima-3}')$ 
may be rewritten in a similar way.

\textbf{Step 3} (tightness). 
Use now Chebyshev inequality and (\ref{stima-2}$'$), \eqref{reg-spaz}.
Then, given $\alpha\in\left(
0,\frac 12\right)  $, $T>0$, for every
$\varepsilon>0$ there is a bounded set $B_{\varepsilon}$ such that
\[
B_{\varepsilon} \subset L^{2}\left(  0,T;V\right)  
                        \cap 
                        W^{\alpha,1+\frac 1b}(0,T;\mathbb W^{-2,1+\frac 1b})
\]
and
\[
 \inf_n P_{n}\left(  B_{\varepsilon}\right)  >1-\varepsilon
\]
The space
$L^2(0,T;V) \cap  W^{\alpha,1+\frac 1b}(0,T;\mathbb W^{-2,1+\frac 1b})$ 
is compactly embedded in $ L^{1+\frac 1b}(0,T;H)$
(see, e.g., Theorem 2.1 in \cite{fg}). 
Hence, for every $\varepsilon>0$
 there is a compact set $K_{\varepsilon}$ such that 
\[
K_{\varepsilon}\subset L^{1+\frac 1b}\left(  0,T;H\right)
 \quad
 \mbox{ and  }
 \quad
 \inf_n P_{n}\left(  K_{\varepsilon}\right)  >1-\varepsilon
\]
Now, take any separable Banach space $\mathbb W^\prime$ such that 
$\mathbb W^{-2,1+\frac 1b}$ is compactly embedded in
$\mathbb W^\prime$;
e.g. $\mathbb W^\prime= \mathbb W^{-2-\theta,1+\frac 1b}$
for some $\theta>0$.
Notice that 
all the spaces 
$W^{1,1 + \frac 1b}(0,T)$, $W^{1,2}(0,T)$, $W^{\alpha,p}(0,T)$ 
(for $\alpha p>1$)
are continuously embedded into $C([0,T])$. Hence, 
the space of vectors with the regularity specified by \eqref{3+reg}
is compactly embedded in $C([0,T];\mathbb W^\prime)$
(see, e.g., Theorem 2.2 in \cite{fg}).
From the boundedness in the mean of $I^n$ 
in $W^{1,1+\frac 1b}(0,T;\mathbb W^{-2,1+\frac 1b})$, of 
$J^{n}$ in $W^{1,2}\left(  0,T;V^\prime \right)  $ 
and of the law of the Wiener process in
$W^{\alpha,p}\left(  0,T;H\right)  $ for every $\alpha\in\left(
\frac 1p,\frac 12\right)  $, again by Chebyshev inequality and compact
embedding 
 we obtain that for every $\varepsilon >0$ 
there exists a compact set $K_{\varepsilon}^{\prime}$ such that
\[
  K_{\varepsilon}^{\prime}
\subset 
  C\left(\left[0,T\right];\mathbb W^{\prime}\right)
 \quad
  \mbox{ and }
 \quad
 \inf_n P_{n}\left(  K_{\varepsilon}^{\prime}\right)  >1-\varepsilon
\]

Therefore the family of measures $\left\{  P_{n}\right\}  $ is
tight in $L^{1+\frac 1b}\left(  0,T;H\right)  $
and in $C\left(  \left[  0,T\right];\mathbb W^{\prime}\right)  $, 
with their Borel $\sigma$-fields.
 Hence there exists
a probability measure $P$ on
\[
 C\left(\left[0,T\right];\mathbb W^{\prime}\right)
    \cap  
 L^{1+\frac 1b}\left(  0,T;H\right)
\]
that is the weak limit in such spaces of a subsequence $\left\{  P_{n_{k}
}\right\}  $.

\textbf{Step 4} ($P$ is a solution to the martingale problem).
From the uniform estimates on $\left\{  P_{n_{k}}\right\}  $ in $L^{2}\left(
0,T;V\right)  $, $L^{\infty}\left(  0,T;H\right)  $ and
$L^{1+b}(0,T;\mathbb L^{1+b})$
 we may deduce that $P$
gives probability one to each one of these spaces and has bounds in the mean
similar to those uniform of $P_{n_{k}}$. 
 This way we have checked   property [MP1] in the definition of
solution to the martingale problem.

Concerning [MP3], we have
$P_{n_{k}}\rightarrow P$ as weak convergence of probability measures on 
$C([0,T];\mathbb W^\prime)$; in particular 
$\Pi_{0}P_{n_{k}}\rightarrow\Pi_{0}P$
as probability measures on $\mathbb W^{\prime}$. 
But $\Pi_{0}P_{n_{k}}$ is
the law of $\pi_{n_k}u_{0}$, which converges to $\mu$ since $\pi_{n_k}u_{0}$
converges $Q$-a.s. to $u_{0}$. Hence $\Pi_{0}P$ is $\mu$.

Finally, let us check property [MP2]. We proceed as in \cite{dpz}
(Sec. 8.4) or in \cite{fg}. 
\\
Given $\psi\in\mathcal{D}^{\infty}$,
we have to prove that for every $t>s\geq0$ and every 
bounded $F_{s}$-measurable
random variable $Z$, we have
\begin{align*}
 E^P \left[  \left(  M_{t}^{\psi}\right)  ^{2}\right]   
&  <\infty
\\
 E^P \left[  \left(  M_{t}^{\psi}-M_{s}^{\psi}\right)  Z\right]   
& =0
\\
 E^P \left[  \left(  \big[(M_t^\psi)^2-\varsigma_t\big]
                    -\big[(M_s^\psi)^2-\varsigma_s\big]\right)Z\right]
&  =0
\end{align*}
where $\varsigma_t := \int_0^t |G(\xi_s)\psi|^2_H ds$.
Defined
\begin{multline*}
 M_t^{\psi,n_k}
  :=\langle \xi_t,\pi_{n_k}\psi\rangle_H-\langle \xi_0,\pi_{n_k}\psi\rangle_H
   -\int_0^t \langle \Phi(\xi_s),\pi_{n_k}A\psi\rangle ds
\\   +\int_0^t \langle B(\xi_s,\pi_{n_k}\psi),\xi_s
         \rangle ds\,,
\end{multline*}
for the measure $P_{n_{k}}$ we
know (see, e.g., \cite{dpz} Sec 8.4 ) that $\left(  M_{t}^{\psi,n_{k}
},F_{t},P_{n_{k}}\right)  $ is a square integrable martingale with quadratic
variation
\[
\left[  M^{\psi,n_{k}}\right]  _{t}
 \equiv \varsigma_t^{n_k}
 =
 \int_0^t |\pi_{n_k}G(\xi_s)\psi|^2_H ds
\]
Thus
\begin{align}
 E^{P_{n_k}} \left[  \left(  M_{t}^{\psi,n_{k}}-M_{s}^{\psi,n_{k}}\right)
  Z\right]   &  =0 \label{m1}
\\
 E^{P_{n_k}} \left[\left( \big[ (M_t^{\psi,n_k})^2-\varsigma_t^{n_k}\big]
      -\big[ (M_s^{\psi,n_k})^2-\varsigma_s^{n_k}\big]  \right) Z\right]   
 &  =0\label{m2}
\end{align}
Moreover, by \eqref{2e-marting} we know that there exists some
$\varepsilon>0$  such that
\begin{equation}\label{2+e}
 \sup_k E^{P_{n_k}}\left|M_t^{\psi,n_k}\right|^{2+\varepsilon}<\infty
\end{equation}

Now, let us consider the limit as $k \to \infty$.
\\
We know that $P_{n_k}$ converges weakly to $P$; then by Skorohod
theorem
there exists a stochastic basis $(\tilde\Omega, \tilde F,\tilde F_t,
\tilde P)$ and, on this basis, there exist
$L^{1+\frac 1b}(0,T;H)\cap C([0,T];\mathbb W^\prime)$-valued random 
variables $\tilde u, \tilde u^{n_k}$ such that 
$\tilde u$ has the same law of $u$, $\tilde u^{n_k}$ has the same
law of
$u^{n_k}$ and $\tilde u^{n_k} \to \tilde u$ $\tilde P$-a.s. in the 
$L^{1+\frac 1b}(0,T;H)\cap C([0,T];\mathbb W^\prime)$-norm.

Define 
\begin{multline*}
 \tilde  M_t^{\psi,n_k}
  :=\langle \tilde u_t^{n_k},\psi\rangle_H
   -\langle\tilde u_0^{n_k},\psi\rangle_H
   -\int_0^t \langle \Phi(\tilde u_s^{n_k}),A\psi\rangle ds
\\   +\int_0^t \langle B(\tilde u_s^{n_k},\pi_{n_k}\psi),\tilde u_s^{n_k}
         \rangle ds
\end{multline*}
Then \eqref{m1}-\eqref{2+e} hold true (with the obvious change of notation).
\\
If we prove that $\tilde M_t^{\psi,n_k}\to \tilde M_t^{\psi}$ $\tilde
P$-a.s. as $k \to \infty$,
then by the equiboundedness relationship \eqref{2+e} we obtain that 
$\tilde M_t^{\psi,n_k}\to \tilde M_t^{\psi}$  in $L^1(\tilde \Omega,\tilde P)$ 
and in $L^2(\tilde \Omega,\tilde P)$ and $\tilde \varsigma_t^{n_k}
\to \tilde \varsigma_t$ in 
$L^1(\tilde \Omega,\tilde P)$.
This concludes our proof.
So, we have to prove a $\tilde P$-a.s. convergence for each term in the
definition
of $\tilde M_t^{\psi,n_k}$.

It is trivial that $\tilde P$-a.s.
$$
 \langle \tilde u_t^{n_k},\psi\rangle_H-\langle
 \tilde u_0^{n_k},\psi\rangle_H
 \;\to\;
 \langle \tilde u_t,\psi\rangle_H-\langle \tilde u_0,\psi\rangle_H 
$$
Notice that there appears the scalar product in $H$
and not the duality pairing $\langle \tilde u_t,\psi\rangle_{\mathbb
  W^\prime,\mathbb W}$, because 
the limit process $u$  belongs
to $C([0,T];H_\sigma)$ with probability one.

Moreover, there exists a subsequence (we do not write that we consider
a subsequence, since we shall pass through subsequences a few times
from now on)
such that
$$
 \tilde P-\text{a.s.} \qquad
 \tilde u_s^{n_k}(x) \to \tilde u_s(x) \quad \text{ for a.e. } (s,x) \in
 [0,T]\times \mathcal T 
$$
We also have, for any $k$
$$
 \tilde E \int_0^T |\tilde u^{n_k}_s|^{1+b}_{\mathbb L^{1+b}}ds\le C_3
$$
Keeping in mind \eqref{1+1b} and \eqref{BconL4}, it follows that $\Phi(\tilde
u^{n_k})$
is equibounded in $L^{1+\frac 1b}(\tilde \Omega\times [0,T]\times
\mathcal T)$ and
$\langle B(\tilde u^{n_k},\pi_{n_k}\psi),\tilde  u^{n_k}\rangle$ is 
equibounded in $L^{\frac{1+b}2}(\tilde \Omega\times [0,T])$ respectively.
First, we get that $\Phi(\tilde u^{n_k}_s(x))
\to \Phi(\tilde u_s(x))$ $\tilde P$-a.s. and for a.e. $(s,x)$.
By the equiboundedness of $\Phi(\tilde u^{n_k}_s(x))$ in 
$L^{1+\frac 1b}(\tilde \Omega\times [0,T]\times \mathcal T)$
it follows that $\Phi(\tilde u^{n_k})$ converges to $\Phi(\tilde u)$
in $L^1(\tilde \Omega\times [0,T]\times\mathcal T)$; 
we get
$$
 \tilde E \int_0^t \langle \Phi(\tilde u^{n_k}_s),A\psi\rangle ds
 \to
 \tilde E \int_0^t \langle \Phi(\tilde u_s),A\psi\rangle ds
$$
Hence a subsequence of 
$\int_0^t \langle \Phi(\tilde u^{n_k}_s),A\psi\rangle ds$ 
converges $\tilde P$-a.s.

On the other hand, another (sub)subsequence can be extracted so that
$$
 \tilde P-\text{a.s.}\qquad 
 \tilde u_s^{n_k} \to \tilde u_s \quad 
  \text{ in } \mathbb L^4 \text{ for a.e. } s 
$$
Then,  by triangle inequality
$\langle B(\tilde u^{n_k}_s,\pi_{n_k}\psi),\tilde  u_s^{n_k}\rangle
\to \langle B(\tilde u_s,\psi),\tilde  u_s\rangle$ $\tilde P$-a.s. and for
a.e. $s$.
By the equiboundedness 
of $\langle B(\tilde u^{n_k}_s,\pi_{n_k}\psi),\tilde  u_s^{n_k}\rangle$ in 
$L^{\frac {1+b}2}(\tilde \Omega\times [0,T])$, 
we conclude as above that there exists a subsequence 
of $\int_0^t \langle B(\tilde u^{n_k}_s,\psi),\tilde  u^{n_k}_s\rangle ds$
converging $\tilde P$-a.s. .

Considering the convergence of a suitable subsequence (the last
extracted), we get that \eqref{m1}-\eqref{2+e} 
in the limit allow to conclude the proof.

\textbf{Step 5} (uniqueness). Let $u^{(i)}$, $i=1,2$, be two
strong  solutions on the same Brownian stochastic basis. 
We are going to prove pathwise uniqueness, which implies uniqueness of
martingale solutions.

Let
\begin{align*}
v_t  & =u^{(1)}_t - u^{(2)}_t,
          \quad v^m_t=\pi_m v_t\\
\theta_t    
      & =2 C_B \left[|u^{(1)}_t|_{\mathbb L^5}^5
         + |u^{(2)}_t|_{\mathbb L^5}^5\right]
        +L_G.
\end{align*}
with $C_B$ as in Lemma \ref{1} and $L_G $ as in \eqref{Glip}.

We have
\begin{equation}
\label{dVprime}
  de^{-\int_0^t \theta_s ds}\left|v^m_t\right|
   _{V^\prime}^2
\\
  =-\theta_t  e^{-\int_0^t \theta_s ds}
              \left|v^m_t\right|_{V^\prime}^2 dt
   +e^{-\int_0^t \theta_s  ds}
    d\left|v^m_t\right|_{V^\prime}^2
\end{equation}
By It\^{o} formula, 
the last differential is
\begin{align*}
 d\left|v^m\right|_{V^\prime}^2  
= &
 -2\langle 
 \pi_m [\Phi(u^{(1)})-\Phi(u^{(2)})],v^m\rangle _H dt
 - 
 2\langle \pi_m[ B(u^{(1)},v)+B(v,u^{(2)})],A^{-1}v^m \rangle_H dt
\\
 & +2\langle \pi_m [G(u^{(1)})-G(u^{(2)})] \,dw,A^{-1}v^m\rangle_H
 +
 \left\| \pi_m A^{-1/2} [G(u^{(1)})-G(u^{(2)})] \right\|_{HS(H)}^{2}dt
\end{align*}
We estimate some terms as follows.
By Lemma \ref{1}
\begin{align*}
2|\langle \pi_m[ B(u^{(1)},v)+B(v,u^{(2)})],A^{-1}v^m \rangle|
&=
2|\langle B(u^{(1)},v)+B(v,u^{(2)}),A^{-1}v^m \rangle |
\\
&\le
\nu \left| v\right|_H^2 + 2C_B\left( |u^{(1)}_t|_{\mathbb L^5}^5
         + |u^{(2)}_t|_{\mathbb L^5}^5\right) 
        \left|v^m\right|^2_{V^\prime}
\end{align*}
By Lemma \ref{3}
\begin{align*}
2\langle  \pi_m [\Phi(u^{(1)})-\Phi(u^{(2)})],v^m\rangle _H 
&=
2\langle \Phi(\pi_m u^{(1)})-\Phi(\pi_m u^{(2)}),v^m\rangle _H 
-2 \epsilon^m
\\
&
\ge 
 2 \nu \left|v^m \right|_H^2 -2 \epsilon^m
\end{align*}
where $\int_0^T|\epsilon^m_t| dt \le 
      C(\left|u^{(1)}\right|_{L^5(0,T;\mathbb L^5)}, 
        \left|u^{(2)}\right|_{L^5(0,T;\mathbb L^5)})$
and ${\displaystyle\lim_{m \to \infty}} \int_0^T \epsilon^m_t \,dt=0$.
\\
By \eqref{Glip}
\begin{align*}
 \left\| \pi_m A^{-1/2} [G(u^{(1)})-G(u^{(2)})] \right\|_{HS(H)}^2
& \le
 \left\| A^{-1/2} [G(u^{(1)})-G(u^{(2)})] \right\|_{HS(H)}^2
\\
& \le
  L_G  \left|v\right|^2_{V^\prime}
\end{align*}

Now, we integrate in time equation \eqref{dVprime} and use
the above estimates, obtaining
\begin{align*}
& 
  e^{-\int_0^T \theta_s ds}\left|v^m _T\right|_{V^\prime}^2
 +\nu \int_0^T e^{-\int_0^t \theta_s  ds}
      (2\left|v^m_t\right|_H^2-\left|v_t\right|_H^2) dt
\\
& \;
 \leq
 \left|v^m_0\right|_{V^\prime}^2
 +
 2 \int_0^T e^{-\int_0^t \theta_s ds} \epsilon^m_t\,dt
 +
 \int_0^T L_G e^{-\int_0^t \theta_s ds}
          ( \left|v_t\right|_{V^\prime}^2-
             \left|v^m_t\right|_{V^\prime}^2) dt
\\
&\qquad
 +
  2\int_0^T 
   \langle \pi_m [G(u^{(1)}_t)-G(u^{(2)}_t)] \,dw_t,A^{-1}v^m_t\rangle_H 
\end{align*}
We can take the limit as $m \to \infty$ in every term. We get
\begin{align*} 
  e^{-\int_0^T \theta_s  ds}\left|v_T\right|_{V^\prime}^2
& +\nu \int_0^T e^{-\int_0^t \theta_s ds}
      \left|v_t\right|_H^2 dt
\\
&
 \leq
 \left|v_0\right|_{V^\prime}^2
  +
  2\int_0^T \langle [G(u^{(1)}_t)-G(u^{(2)}_t)] \,dw_t,A^{-1}v_t \rangle_H 
\end{align*}
Hence
$$
 E \left[e^{-\int_0^T \theta_s  ds}\left|v_T\right|_{V^\prime}^2\right]
+
  \nu E \left[\int_0^T e^{-\int_0^t \theta_s  ds}
      \left|v_t\right|_H^2 dt \right]
 \leq
 E \left|v_0 \right|_{V^\prime}^2
$$
When the initial conditions of $u^{(i)}$ coincide, we deduce
$$
\int_0^T e^{-\int_0^t \theta_s  ds}
       \left|  v_t\right|_H^2 dt=0
$$
with probability one. Since $\int_{0}^{T}\theta_s  ds<\infty$
a.s., we have $v=0$ a.s., as considering $v$ as a measurable function of $t$
with values in $H$. This implies
that with probability one $u^{(1)} =u^{(2)}$,
where the equality holds in $L^\infty(0,T;H)$.  
\end{proof}

\subsection{Markov and Feller property}

\begin{lemma}\label{mark}
Let $u_0^n, u_0$ be initial data satisfying the assumption of Theorem 
\ref{teoexist1} and let $(u_t^n)_{t\geq 0}$ and $(u_t)_{t\geq0}$ be 
the corresponding strong solutions on the
same given Brownian stochastic basis. 
\\
If $E|u_0^n- u_0|^2_{V^\prime} \rightarrow 0$, then 
for every $T>0$, $(u_t^n)_{t\geq 0}$ converges to $(u_t)_{t\geq 0}$ in
probability on $\left[  0,T\right]  \times\Omega$ in the topology of
$H$, and $u_{T}^{n}$ converges to $u_{T}$
in probability on $\Omega$ in the topology of $V^\prime$.
\end{lemma}

\begin{proof}
We proceed as in the previous Step 5 to get the following estimate
$$
 E \left[e^{-\int_0^T \theta_s  ds}
  \left|u_T^n-u_T \right|_{V^\prime}^2\right]
+
  \nu E \left[\int_0^T e^{-\int_0^t \theta_s ds}
      \left| u_t^n-u_t \right|_H^2 dt \right]
 \leq
 E \left| u_0^n-u_0 \right|_{V^\prime}^2
$$
where
$$
 \theta_t
=
 2 C_B \left[|u^n_t\right|_{\mathbb L^5}^5
         + \left|u_t|_{\mathbb L^5}^5\right] + L_G
$$
Since $\int_{0}^{T}\theta_s  ds<\infty$ with probability one, we
get the result.
\end{proof}

\begin{theorem}
The strong solutions of equation (\ref{abstr SNS}) on a given Brownian stochastic
basis define a Markov process in $H$ with the Feller property in $V^\prime$.
\end{theorem}
\begin{proof}
Denote by $u(t;y)$ the solution at time $t$ which started at time 0 
from $y$.

Given $t>0$ the dynamics $y \mapsto u(t;y)$ is uniquely defined in
$H$; hence the 
Markov property is inherited by $u$ from the Galerkin approximations
$u^n$.

The process solution enjoys the Feller property if
$$
 E g(u(t;z)) \to E g(u(t;y)) \qquad \text{as } z \to y \text{ in }V^\prime
$$
for any $t\ge 0, g \in C_b(V^\prime)$.
For this it is enough the convergence in probability: 
$u(t;z) \to u(t;y)$ as $z \to y $ in $V^\prime$.
But, as in Lemma \ref{mark} (now the initial data are deterministic), 
we know that 
$$
 E\left[  e^{-\int_0^t \theta_s ds}|u(t;z)-u(t;y)|_{V^\prime}^2\right]  
\le 
  |z-y|^2_{V^\prime}
$$
Then, we conclude as before that $|u(t;z)-u(t;y)|_{V^\prime} \to 0$ in
probability as $|z-y|_{V^\prime} \to 0$.
\end{proof}

\subsection{Stationary solutions}

As in \cite{fg}, existence of stationary solutions is obtained in the
limit, showing first that the Galerkin problem has at least one 
stationary solution.
Our result is the following

\begin{theorem} \label{teostaz}
Assume that $2\nu \frac{(2\pi)^2}{L^2}> \lambda_0$.
Then equation \eqref{abstr SNS} has a stationary solution.
\end{theorem}
\begin{proof}
Let us consider
$$
  du_t^{n}+[A\Phi(u_t^{n})+\pi_{n}B(u_t^{n},u_t^{n})]\,dt
=
 \pi_{n}G\left(u_t^{n}\right)  \,dw_t,
     \quad u^{n}_0=0
$$
We use estimates from Appendix 1.
By \eqref{prima stima}, using $\|u\|_V \ge \frac{2\pi}L |u|_H$
we get
$$
\begin{array}{rl}
  \displaystyle
 \frac{d\;}{dt} E |u^n_t|^p_H +p\nu \tfrac{(2\pi)^2}{L^2} E|u^n_t|^p_H
 &\le
  \frac 12 p(p-1) \left[\lambda_0 E |u^n_t|^p_H + \rho E |u^n_t|^{p-2}_H\right]
\\[2mm]
 &\le  
  \frac 12 p(p-1)(\lambda_0+\varepsilon)E |u^n_t|^p_H +C(\varepsilon,p,\rho)
\end{array}
$$
for some positive $\varepsilon$.
\\
If $2\nu \tfrac{(2\pi)^2}{L^2}> \lambda_0$, then 
there exist $p>2$ and $\varepsilon>0$ such that $p\nu  \tfrac{(2\pi)^2}{L^2}
>\frac 12 p(p-1)(\lambda_0+\varepsilon)$.
Therefore there exists $a>0$ such that
$$
  \frac{d\;}{dt} E |u^n_t|^p_H + a  E |u^n_t|^p_H
  \le C(\varepsilon,p,\rho) 
  \qquad \text{ with } u^n_0=0 ;
$$
by Gronwall Lemma we get
$$
 E|u^n_t|_H^p \le C_6 \qquad \forall t\ge 0, \forall n\ge 1
$$
Hence, the  family of random variables $\{u^n_t\}_{t \ge 0}$ is  tight in
$H_n$.
Notice that the Galerkin problem is Feller in $H_n$.
Then, by the Krylov-Bogoliubov method we get that there exists a stationary
solution (whose law we denote by $\mu_n$)
for the Galerkin equation.

Now, consider the Galerkin problem with initial velocity of law
$\mu_n$ and denotes the law of the solution by $P_n$ (a probability
measure on $C([0,\infty);\mathbb W^\prime)$.
We have
$$
 E^{P_n}|\xi_0|_H^p \le C_6 \qquad \forall n \ge 1
$$
The corresponding solution $P_n$ is a stationary process in $H_n$, i.e.
$$
 P_n\left\{ \xi_t +\int_r^t A\Phi(\xi_s)ds + 
   \int_r^t \pi_n B(\xi_s,\xi_s)ds= \xi_r+\int_r^t \pi_n
 G(\xi_s)dw_s\right\}=1
$$
and
$$
 P_n(\xi_t)=P_n(\xi_r)
$$
for any $ 0 \le r \le t < \infty$.

Now we proceed as in \cite{fg}. Endow 
$L^{1+\frac 1b}(0,\infty;H)$ with the distance
$$
 d_{1+\frac 1b}(u,v)=\sum_{k=1}^\infty \frac 1 {2^k}
     \left( |u-v|_{L^{1+\frac 1b}(0,k;H)} \wedge 1\right)    
$$
and $C([0,\infty];\mathbb W^\prime)$ 
with the distance
$$
 d_\infty(u,v)= \sum_{k=1}^\infty \frac 1 {2^k}
     \left( |u-v|_{C([0,k];\mathbb W^\prime)} \wedge 1\right)
$$
The convergence with respect to  $d_{1+\frac 1b}+d_\infty$
is equivalent with the convergence on every finite time interval.
We come back to the bounds $(\ref{mediaP}')$, $(\ref{stima-2}')$, 
$(\ref{stima-3}')$ and \eqref{reg-spaz}, to notice that they hold
true because they depend only on $E |u_0^n|^p$. Thus we get tightness
on every finite interval; we pass to the limit for a subsequence and
get the limit process $P$ which is stationary, since the $P_n$ are so.
It can be shown as  before that $P$ is a martingale solution 
to equation \eqref{abstr SNS}. 

Defined the Markov semigroup $\mathbf P_t$ acting on the space 
of Borel bounded functions $B_b(H)$ as $\mathbf P_t \phi(y)=E\phi(u(t;y))$,
we get
that the law $\mu$ of this stationary solution is an invariant measure, in the
sense that
$
 \int \mathbf P_t \phi \, d\mu =  \int \phi \, d\mu
$
for any  $\phi \in B_b(H)$ and $t \ge 0$.
\end{proof}

\section{The case $\Phi(u)=\nu |u|^4u$} \label{b=5}

Instead of \eqref{Ipo-phi}, let us assume that
$$
 \Phi(u)=\nu |u|^4u
$$
This corresponds to the case $b=5$ with the nonlinearity
acting everywhere. 
The interest in this model will be explained in Subsection 4.2.

We can  analyze this model as done in the previous
section,
with  few changes.
Mainly, the solution will live in 
$L^\infty(0,T;H) \cap L^6(0,T;\mathbb X)$, where 
$\mathbb X$ is the closure of $\mathcal D^\infty$ 
w.r.t. the norm 
$$
 |u|_{\mathbb X}=
  \left(\int_{\mathcal T}\{|u(x)|^4 |\nabla u(x)|^2 +4 |u(x)|^2
  \sum_{i=1}^3 [u(x)\cdot \partial _i u(x)]^2 \}dx
  \right)^{1/6}
$$
Notice that the term 
$ \int_0^t \left\langle\Phi(u_s),A\psi\right\rangle_H ds$
in the equation is well defined, since
$$
 \int_0^t \left|\left\langle\Phi(u_s),A\psi\right\rangle \right| ds
\le |A\psi|_{\mathbb L^6}
   \int_0^t \left|\Phi(u_s)\right|_{\mathbb L^{6/5}} ds
 =
    |A\psi|_{\mathbb L^6} \nu \int_0^t \left|u_s\right|_{\mathbb L^6}^5 ds
$$
The last integral is well defined for functions $u \in L^5(0,T;\mathbb L^6)$.
But, if $u \in \mathbb X$, then $u \in \mathbb L^6$ by Theorem
\ref{compatt}
in Appendix 2.

We have the following result
\begin{theorem} \label{teob5}
Let $\Phi(u)=\nu |u|^4u$.\\
Let $\mu$ be a measure on $H$ such that 
$m_p:=\int_H |v|_H^p \mu(dv)  <\infty$ for some $p>2$. 
Then there
exists at least one solution to the martingale problem \eqref{weakSNSeq}
with initial
condition $\mu$, assuming that  condition $[MP1]$ in Definition 5
is replaced with 
\[
 P\left(  \sup_{t\in [0,T]} \left|\xi_t\right|_H
  +\int_0^T \left\|\xi_s\right\|_{\mathbb X}^6 ds
   <\infty\right)  =1
\]
Moreover, if $2 \nu C_X >\lambda_0$ (with 
$|u|^6_{\mathbb X} \ge C_X |u|^6_H$), then
there exists a stationary solution.
\end{theorem}
\begin{proof}
Let us check step by step how our previous proof (for $b=5$) can be
adapted to handle this model.\\
\textbf{Step 1} Instead of Lemma \eqref{2}, we use  
$$
   \langle A\Phi(u),u\rangle_H = \nu |u|_{\mathbb X}^6
$$
Hence if we  apply It\^{o} formula  (for $p\ge 2$) 
to $|u^n_{t\wedge \tau^n_R} |_H^p$, we get \eqref{mediaP}; however
\eqref{stima-2} is replaced by 
$$
   E \int_0^T \|u^n_{s\wedge\tau_R}\|_{\mathbb X}^6 1_{\{s< \tau_R\}}ds
\le C_2^\prime
$$
\eqref{stima-3} is
a consequence of the latter relationship, since $ \mathbb X
\subset \mathbb L^6$.
\\
\textbf{Step 2} 
The estimates are still valid:
$$
 \sup_n E|u^n|_{W^{\alpha,\frac 65}(0,T;\mathbb W^{-2,\frac 65})}<\infty
   \qquad \text{ for } 0<\alpha <\tfrac 12
$$
and 
$$
 \sup_n E| u^n|_{L^6(0,T;\mathbb X)}<\infty
$$
\textbf{Step 3} 
 What we need is
a compact embedding, which is given in Theorem \ref{compatt}.
Thus, by Theorem 2.1 in \cite{fg} the space
$L^6(0,T;\mathbb X)\cap 
W^{\alpha,\frac 65}(0,T;\mathbb W^{-2,\frac 65})$ is compactly embedded in 
$L^{\frac 65}(0,T;\mathbb L^6)$.
Therefore the family of measures $\{P_n\}$ is tight in 
$L^{\frac 65}(0,T;\mathbb L^6)$ and in $C([0,T];\mathbb W^\prime)$,
chosen $\mathbb W^\prime$ such that $\mathbb W^{-2,\frac 65}$ 
is compactly embedded in $\mathbb W^\prime$.
\\
\textbf{Step 4} 
The remaining part of the proof for the existence holds true.
\\
\textbf{Step 5} 
As far as the uniqueness is concerned, Lemma 
\ref{3} has to be replaced with 
$$
 \langle\Phi(u^{(1)})-\Phi(u^{(2)}),u^{(1)}-u^{(2)}\rangle_H
\ge 
 0
$$
This comes from its proof, when we put $\sigma=\tilde \sigma$.
The above inequality is not enough to get uniqueness.
Other estimates failed to be useful so far and uniqueness is an open problem.
\\
\textbf{Stationary martingale solutions}.
We consider the sequence of Galerkin solutions $\{u^n\}_{n \ge 1}$, all with
zero initial velocity.
From the estimates in the Appendix 1, we get
$$
  \frac{d\;}{dt} E |u^n_t|^p_H +p\nu  E[|u^n_t|^{p-2}_H |u^n_t|^6_{\mathbb X}]
 \le
  \tfrac 12 p(p-1)\lambda_0 E |u^n_t|^p
  + 
  \tfrac 12 p(p-1) \rho E|u^n_t|^{p-2}_H 
$$
By the embeddings $\mathbb X \subset \mathbb L^6\cap H $,
we get $|u|^6_{\mathbb X} \ge C_X |u|^6_{H}$; thus
$$
  \frac{d\;}{dt} E |u^n_t|^6_H + p\nu C_X E [|u^n_t|^{p-2}_H |u^n_t|^6_H]
  \le
  \tfrac 12 p(p-1)\lambda_0 E |u^n_t|^p
  + 
  \tfrac 12 p(p-1) \rho E|u^n_t|^{p-2}_H 
$$  
Using that $|u|^2_H \le |u|^6_H + \frac 2{3\sqrt 3}$, we obtain by
easy computations that 
$$
 \frac{d\;}{dt} E |u^n_t|^6_H + p\nu C_X E |u^n_t|^p_H 
  \le
  \tfrac 12 p(p-1)(\lambda_0+\varepsilon) E |u^n_t|^p
  + C(\varepsilon, p, \rho, C_X, \nu)
$$
for any $\varepsilon >0$ (with the latter constant $C$ being a suitable positive
constant).
If $2 \nu C_X > \lambda_0$, we conclude as before 
by Gronwall Lemma that there exists $p>2$ such that 
$$
  E |u^n_t|^p_H \le C_7 \qquad
   \forall t\ge 0, \forall n \ge 1
$$
From now on, the proof proceeds as in the previous case.
\end{proof}

\subsection{Scaling}

Let us start with an heuristic digression.
We recall that in the Kolmogorov 1941 theory of turbulence
(see, e.g., Sect. 6.3.1 in \cite{f}  dealing with the deterministic
equations
and \cite{k} dealing with the stochastic equations), 
one believes
that the following equality in law is approximatively true
$$
 u(r+\lambda x)-u(r) \;\stackrel{\text{in law}}{=} \; 
 \lambda^{1/3} [u(r+x)-u(r)]
$$
for any $r,x \in \mathbb R^3$ and for $\lambda$ in some range of small
positive real numbers.
\\
(In the whole section, $u(x)$, without the time variable,
 denotes a stationary solution.) 
\\
This implies
$ \lambda^{-1/3} [u(\lambda x)-u(0)] \;\stackrel{\text{in law}}{=} \; 
   u(x)-u(0)$.
\\
According to this result, we are interested in the scaled velocity
$u_\lambda$, defined
by the following scaling transformation
\[
u_{\lambda}(t,x):=\lambda^{-1/3}u\left(  \lambda^{2/3}t,\lambda x\right)
\]
for $\lambda\in\left(  0,1\right)  $; hence the function
$u_{\lambda}(t,x)$ is defined for 
$x\in\left[
0,\frac{L}{\lambda}\right]  ^{3}$.

We assume that $u=u(t,x)$ solves in the torus $[0,L]^3$ 
the modified Navier--Stokes equation 
(with $\Phi(u)=\nu |u|^4u$),
with additive noise 
$$
 du+\left[  -\triangle\Phi(u) + (u\cdot\nabla)u
   +\nabla q\right]  dt   =\sum_{(k,j) \in \Lambda} 
  \sigma_{k,j}  d\beta_{k,j} (t)  h_{k,j}
$$
(with respect to \eqref{cilindrico}, there are the coefficients
$\sigma_{k,j}$; some of them may vanish and
therefore we denote by $\Lambda$ the set for the summation on 
$\sigma_{k,j}\neq 0$. Condition \eqref{crescG} is satisfied if
$\sum_{(k,j)\in\Lambda} |\sigma_{k,j}|^2 <\infty$.)
\\
The scaled velocity  satisfies an equation very similar to this one.
\begin{proposition}
We have
$$
du_{\lambda}+\left[  -\triangle\Phi\left(  u_{\lambda}\right)
+\left(  u_{\lambda}\cdot\nabla\right)  u_{\lambda}+\nabla q_{\lambda}\right]
dt
=
  \sum_{(k,j)\in\Lambda}
  \sigma_{k,j}  d\beta^\lambda_{k,j} (t)  h^\lambda_{k,j}$$
where  $q_{\lambda}$ is a suitable function,
$h^\lambda_{k,j}(x)=h_{k,j}(\lambda x)$ and the processes 
\[
\beta_{k,j}^{\lambda}\left(  t\right)  :=\lambda^{-1/3}\beta_{k,j}\left(
\lambda^{2/3}t\right)
\]
are independent standard Brownian motions.
\end{proposition}

\begin{proof}
The rigorous
proof has to be performed at the level of the integral weak
formulation  of the equation
and it is tedious and elementary. We just point out the main (somewhat
heuristic) arguments behind it.
We have
\[
\left.  \frac{\partial u_{\lambda}}{\partial t}\right|  _{(t,x)}=\lambda
^{1/3}\left.  \frac{\partial u}{\partial t}\right|  _{\left(  \lambda
^{2/3}t,\lambda x\right)  }
\]
\[
\left.  \left(  u_{\lambda}\cdot\nabla\right)  u_{\lambda}\right|
_{(t,x)}=\lambda^{1/3}\left.  \left(  u\cdot\nabla\right)  u\right|  _{\left(
\lambda^{2/3}t,\lambda x\right)  }
\]
\begin{align*}
\left.  \frac{\partial\beta_{k,j}^{\lambda}}{\partial t}\right|  _{(t)}  &
\sim\frac{\beta_{k,j}^{\lambda}\left(  t+dt\right)  -\beta_{k,j}^{\lambda}\left(
t\right)  }{dt}\\
& =\lambda^{-1/3}\frac{\beta_{k,j}^{\lambda}\left(  \lambda^{2/3}t+\lambda
^{2/3}dt\right)  -\beta_{k,j}^{\lambda}\left(  \lambda^{2/3}t\right)  }{dt}\\
& =\lambda^{1/3}\frac{\beta_{k,j}\left(  \lambda^{2/3}t+\lambda^{2/3}dt\right)
-\beta_{k,j}\left(  \lambda^{2/3}t\right)  }{\lambda^{2/3}dt}\\
& \sim\lambda^{1/3}\left.  \frac{\partial\beta_{k,j}}{\partial t}\right|
_{(\lambda^{2/3}t)}
\end{align*}
$$
 \left. \triangle\Phi_{\lambda}\left(  u_{\lambda}\right)  \right|  _{(t,x)}
 =\lambda^{1/3}\left.  \triangle\Phi\left(  u\right)  \right|  _{\left(
\lambda^{2/3}t,\lambda x\right)  }
$$
because
\begin{align*}
D_{i}\left|  u_{\lambda}\right|  ^{4}u_{\lambda}  & =\lambda^{-5/3}
D_{i}\left[  \left|  u\left(  \lambda^{2/3}t,\lambda x\right)  \right|
^{2\cdot2}u\left(  \lambda^{2/3}t,\lambda x\right)  \right]  \\
& =\lambda^{-5/3}u\left(  \lambda^{2/3}t,\lambda x\right)  2\left|  u\left(
\lambda^{2/3}t,\lambda x\right)  \right|  ^{2}u\left(  \lambda^{2/3}t,\lambda
x\right)  \cdot D_{i}\left[  u\left(  \lambda^{2/3}t,\lambda x\right)
\right]  \\
& +\lambda^{-5/3}\left|  u\left(  \lambda^{2/3}t,\lambda x\right)  \right|
^{2\cdot2}D_{i}\left[  u\left(  \lambda^{2/3}t,\lambda x\right)  \right]  \\
& =\lambda\cdot\lambda^{-5/3}\left.  \left[  u2\left|  u\right|  ^{2}
u\cdot\left(  D_{i}u\right)  +\left|  u\right|  ^{2\cdot2}\left(
D_{i}u\right)  \right]  \right|  _{\left(  \lambda^{2/3}t,\lambda x\right)
}\\
& =\lambda\cdot\lambda^{-5/3}\left.  D_{i}\left[  \left|  u\right|
^{4}u\right]  \right|  _{\left(  \lambda^{2/3}t,\lambda x\right)  }
\end{align*}
and then
\[
\left.  \triangle\left|  u_{\lambda}\right|  ^{4}u_{\lambda}\right|
_{(t,x)}=\lambda^{2}\cdot\lambda^{-5/3}\left.  \triangle\left[  \left|
u\right|  ^{4}u\right]  \right|  _{\left(  \lambda^{2/3}t,\lambda x\right)  }
\]
\end{proof}

Now, for any $r \in \mathbb R^3$, 
define the (space) translation operator $\hat r$ as $(\hat r V)(x)=V(x+r)$, to
be understood as 
an identity in the distributional sense. 
We say that a process $V$ is spatially homogeneous if all the space
increments $\delta V(x,h)=V(x+h)-V(x)$ are
statistically invariant with respect to the translation operator $\hat r$:
$V(x+h)-V(x) \;\stackrel{\text{in law}}{=} \; V(x+h+r)-V(x+r)$.
\\
In the same way, we say that a process $V$ is isotropic 
if the law of all the space increments $\delta V(x,h)$ do
not change under simultaneous rotation $\theta$ of the space variables and of
the vector $V$.
Since the space variable lives in a torus, only
rotations of $\mathbb{R}^{3}$ which leave the torus invariant are allowed.

It is easy to check when the  Wiener process on the r.h.s. of our equation
enjoys these statistical invariances. 
Indeed 
$$
 E|w(t,x+h)-w(t,x)|^2=
 t\sum_{(k,j)\in\Lambda}|\sigma_{k,j}|^2 
 |h_{k,j}(x+h)-h_{k,j}(x)|^2
$$
If for any $k$ with $(k,j)\in \Lambda$, 
each coefficient $\sigma_{k,j}$ in front of 
$\cos(\frac{2\pi}L k \cdot x)$ is
equal, in absolute value, to a coefficient in front of 
$\sin(\frac{2\pi}L k \cdot x)$, then 
the second moment of the space increment  is equal to  
$t\frac 2{L^3}\displaystyle \sum_{\substack{(k,j)\in\Lambda\\j=1,2}}
|\sigma_{k,j}|^2 |e^{i\frac{2\pi}{L} k\cdot h}-1|^2$
and therefore depends only on $h$. This implies that $w$ is 
spatially homogeneous (considering the space
variables in $\mathbb R^3/[0,L]^3$). 
On the other hand,
$$
 E|\theta w(t,\theta x+\theta h)-\theta w(t,\theta x)|^2=
 t\sum_{(k,j)\in\Lambda}|\sigma_{k,j}|^2 
 |h_{k,j}(\theta x+\theta h)-h_{k,j}(\theta x)|^2
$$ 
But $h_{k,j}(\theta x+\theta h)-h_{k,j}(\theta x)=
h_{\theta^{-1}k,j}(x+h)-h_{\theta^{-1}k,j}(x)$. 
Then we can
consider only rotations $\theta$ such that 
$(k,j) \in \Lambda \iff (\theta k,j) \in \Lambda$ and in these cases 
$w$ is isotropic  
if $|\sigma_{\theta k,j}|=|\sigma_{k,j}|$ for all $(k,j) \in \Lambda$.

\begin{corollary}
Let $\Lambda$  be such that the process 
$\sum_{(k,j)\in\Lambda}\sigma_{k,j} \beta_{k,j}(t)  h_{k,j}(x)$  is spatially
homogeneous and isotropic. 
For $\Phi(u)=\nu |u|^4u$, 
consider the equation
\begin{equation}\label{NSmod}
\begin{array}{lr} 
 du+\left[  -\triangle\Phi(u) + (u\cdot\nabla)u
 + \nabla q\right]  dt  & =\displaystyle
 \sum_{(k,j)\in\Lambda}
  \sigma_{k,j}  d\beta_{k,j} (t)  h_{k,j}(x) \\[1mm]
 & x  \in\left[  0,L \right]^3
\end{array}
\end{equation}
with initial velocity spatially homogeneous and isotropic,
satisfying the assumptions of Theorem \ref{teob5} and $\sum_{(k,j)\in\Lambda} 
 |\sigma_{k,j}|^2 <\infty$.
Then there exists a solution $u$ spatially homogeneous and isotropic 
for any  $t\geq0$.
For any $\psi \in \mathcal D^\infty$, we have 
\begin{equation}\label{funstr}
E\big[  |\langle u(t,\lambda e) - u(t,0), \psi\rangle |^{p}\big]  
 =
 \lambda^{p/3}E\big[ |\langle u_{\lambda}(\lambda^{-2/3}t,e)
            -u_{\lambda}(\lambda^{-2/3}t,0),\psi \rangle|  ^{p}\big]
\end{equation}
where $u_{\lambda}(t,x)$ is 
spatially homogeneous and isotropic, and solves the
equation 
\begin{equation}\label{NSrisc}
\begin{array}{lr} 
du_{\lambda}+\left[  -\triangle\Phi(u_{\lambda})
  +(u_{\lambda}\cdot\nabla)  u_{\lambda}+\nabla q_{\lambda}\right]
   dt  & =\displaystyle 
 \sum_{(k,j)\in\Lambda}
  \sigma_{k,j}  d\beta^\lambda_{k,j} (t)  h_{k,j}(\lambda x)\\
 & x \in\left[  0,\frac{L}{\lambda}\right]^3
\end{array}
\end{equation}
with initial velocity $u_{\lambda}(0,\cdot)=\lambda^{-1/3}u(0,\lambda \cdot)$.
\end{corollary}

\begin{remark}
The statistical invariance for the solution is obtained from the same
property of the Galerkin approximations, as done  in a similar
context in \cite{vf1}. Indeed, we construct a solution as limit of
a Galerkin subsequence. But the statistical invariance for the Galerkin
processes $u_n$ (for any $n$) is easy to show,
since for any $n$ the finite-dimensional problem has a unique
solution.
\\
Notice that 
in our case we can trivially consider a vanishing  initial
velocity. 
\end{remark}

Now, let us consider the structure function of order $p$
($p=1,2,\dots$)  with respect to a stationary solution
$u$ of the
modified Navier--Stokes equation \eqref{NSmod} :
$$
 S_p(\lambda):=E\left[ |u(\lambda e)-u(0) |^{p}\right] 
$$
($e$ is a unitary vector in $\mathbb R^3$ 
and  $\lambda \in (0,1)$). Similarly, we can work with the
longitudinal structure function.
\\
We point out that Theorem 
\ref{teob5} provides the existence of a stationary solution leaving in
$\mathbb X$,
but this is not enough to define the velocity in {\it every} point of
the torus. We would need to analyze the regularity of stationary
solutions, but we decide to postpone  the study  of existence 
of  more regular stationary solutions to future work
(it would be enough to have the law of $u$  supported
 by the space $C^0(\mathcal T)$).
\\
According to the previous Corollary, 
for the structure function
we get that
\begin{equation}\label{fsL}
 S_p(\lambda)  
 =
 \lambda^{p/3}E\left[|u_\lambda (e)-u_\lambda(0)|^p\right]
\end{equation}
Kolmogorov 1941 theory states (see, e.g.,  \cite{f}, Sect. 6.3.1) that
\begin{equation}\label{K41}
 S_p(\lambda)=C_p \varepsilon^{p/3} \lambda^{p/3}
\end{equation}
where $C_p$ are dimensionless and $\varepsilon$ is the mean energy
dissipation rate.
For $p=2$, \eqref{K41} is the so called two-thirds law of turbulence, 
which is supported by experimental results. 
For $p=3$,  \eqref{K41} is the four-fifths law of turbulence
($C_3=-\frac 45$), deduced from the
assumptions of homogeneity, isotropy and finiteness of the energy
dissipation. For $p>3$,  \eqref{K41}  is not confirmed by
experimental data and its truthfulness is questionable.

According to the above
Proposition and Corollary, we shall provide a relationship similar to
\eqref{K41} for our model \eqref{NSmod}.

Keeping in mind \eqref{fsL}, we investigate the behaviour of 
$E\left[ |u_\lambda(e)-u_\lambda(0)|^p \right]$ 
in order to get insights on the
structure function.
First, we remark that $\beta_{k,j}\left(  t\right)  $ and  $\beta
_{k,j}^{\lambda}\left(  t\right)$ are unitary Brownian motions. 
Thus, the random forces in equations \eqref{NSmod} and \eqref{NSrisc}
are the same in law; what changes is 
the dimension of the torus and correspondingly the eigenvectors
$h_{k,j}$. 
No other terms in the equation (when
projected onto $H$) change with the scaling. It is important to point
out that this property is not true for the usual Navier--Stokes
equation or for the Prouse model introduced at the beginning; namely,
after the scaling the viscous term $\nu \Delta u$ becomes
$\lambda ^{-4/3} \nu \Delta u_\lambda$ (and then the two problems 
for $u$ and $u_\lambda$ are very different, because the 
scaled viscosity $\nu_\lambda=\lambda ^{-4/3} \nu$ 
explodes as $\lambda \to 0$).

Since the coefficients $\sigma_{k,j}$ in the noise 
do not change with  the scaling transformation, 
the mean energy introduced in a unit of volume per unit of time 
by the stochastic forcing
term is independent of $\lambda$ and is equal to 
$\frac 1{L^3}\sum_{(k,j)\in\Lambda}|\sigma_{k,j}|^2$.
 When $\lambda \rightarrow 0$, the size of the domain becomes 
bigger and bigger but the unitary energy does not change.

If, as in turbulence theory, we assume
that during the motion there is  energy transfer from 
large scales to small scales with a universal cascade mechanism depending only
on the unit volume energy, then we would conclude that any  stationary 
state is independent of $\lambda$.
Hence $E\left[ |u_\lambda (e)-u_\lambda(0)|^p \right]
=k_p$ for any $\lambda$.
Coming back to \eqref{funstr}, we conclude that 
$$
 S_p(\lambda) = k_p \lambda^{p/3}
$$
for any $0<\lambda <1$ and $p \in \mathbb N$.

Summing up, we have the following result, providing a result on the
structure function (of any order $p$) under two 
assumptions.
The first assumption is  technical and can be removed as soon as we
are able to prove existence of regular stationary solutions. On the
other hand, the second assumption on energy cascade
has to be considered as an hypothesis quite hard to justify rigorously
(as it is for fluids modeled by the Navier--Stokes equations).
\begin{claim}\label{K41b5}
Let us assume that system \eqref{NSmod} has a stationary solution
$u$, which at any fixed time has  law 
 supported by the space $C^0(\mathcal T)$.
\\
Let us further assume that there is 
energy transfer from 
large scales to small scales with a universal cascade mechanism depending only
on the unit volume energy.

Then for 
 the structure function, given any $0<\lambda <1$ and $p\in \mathbb N$
  we have 
$$
 S_p(\lambda)=k_p \lambda^{p/3}
$$
for some constant $k_p$ independent of $\lambda$.

\end{claim}

\section{Appendix 1: a priori estimates}
We present the estimates on the Galerkin approximations; these are
quite standard (see, e.g., \cite{fg} and \cite{prouse}). Besides the usual
estimates \eqref{mediap}, \eqref{stimaV},
we need also \eqref{stima5} to prove  uniqueness.

Let $(u^n_t)_{t\ge 0}$ be a continuous adapted solution of
equation (\ref{Galerkin SNS}). Let
\[
 \tau_R^n=\inf\left\{  t\ge 0: |u^n_t|_H^2=R\right\}
\]
We have
\begin{align*}
 u^n_{t\wedge \tau^n_R}
&  = 
   u_0^n
  +\int_0^{t\wedge \tau^n_R} \left[
    -A\Phi(u^n_s) -\pi_n B(u^n_s,u^n_s)  \right]  ds
   +\int_0^{t\wedge \tau^n_R} \pi_n G(u^n_s)  dw(s)
\\
& =
  u_0^n
  +\int_0^t \left[
    -A\Phi(u^n_{s\wedge \tau^n_R}) 
    -\pi_n B(u^n_{s\wedge \tau^n_R},u^n_{s\wedge \tau^n_R})  \right] 
       1_{\{s < \tau^n_R\}} ds
\\
&\;
   +\int_0^t 1_{\{s < \tau^n_R\}} \pi_n G(u^n_{s\wedge \tau^n_R})  dw(s)
\end{align*}
For $p\ge 2$ apply It\^{o} formula 
to $|u^n_{t\wedge \tau^n_R} |_H^p$:
\begin{equation*}
\begin{split}
 d|u^n_{t\wedge \tau^n_R}|_H^p 
\le
&
 p |u^n_{t\wedge \tau^n_R}|^{p-2}_H 
      \langle u^n_{t\wedge \tau^n_R}, du^n_{t\wedge \tau^n_R}\rangle_H
\\ 
&
+\frac 12 p (p-1) |u^n_{t\wedge \tau^n_R}|^{p-2}_H
        \|\pi_n G(u^n_{t\wedge \tau^n_R})\|^2_{HS(H)}1_{\{t < \tau^n_R\}}dt
\end{split}
\end{equation*}
Then, integrating in time, we have
\begin{align*}
 | &       u^n_{t\wedge \tau^n_R}|_H^p 
  \le 
 |u_0^n|_H^p
\\
&
+
 p \int_0^t |u^n_{s\wedge \tau^n_R}|_H^{p-2}
    \langle -A\Phi(u^n_{s\wedge \tau^n_R})  
            -\pi_n B(u^n_{s\wedge \tau^n_R},u^n_{s\wedge \tau^n_R})
        ,u^n_{s\wedge \tau^n_R} \rangle_H 1_{\{s <\tau^n_R\}}ds
\\
 &  + p\int_0^t |u^n_{s\wedge \tau^n_R}|_H^{p-2}
 \langle u^n_{s\wedge \tau^n_R}, 
  1_{\{s < \tau^n_R\}}\pi_n G(u^n_{s\wedge \tau^n_R})dw(s) \rangle_H
\\
&
 + \frac 12 p(p-1)\int_0^t |u^n_{s\wedge \tau^n_R}|_H^{p-2} \left\|  \pi_n
   G(u^n_{s\wedge \tau^n_R})  \right\|_{HS(H)}^2 1_{\{s <\tau^n_R\}}ds
\end{align*}
Then, by Lemma \ref{2} and \eqref{incomp}
\begin{align*}
 | u^n_{t\wedge \tau^n_R}|_H^p 
& +
 p\nu \int_0^t |u^n_{s\wedge \tau^n_R}|_H^{p-2} 
              \|u^n_{s\wedge \tau^n_R}\|_V^2 1_{\{s <\tau^n_R\}}ds
\\
&  \le
   \left|  u_0^n\right|_H^p 
  +
  p \int_0^t   |u^n_{s\wedge \tau^n_R}|_H^{p-2}
       \langle u^n_{s\wedge \tau^n_R},1_{\{s < \tau^n_R\}}
       \pi_n G(u^n_{s\wedge \tau^n_R}) dw(s) \rangle_H 
\\
& + 
    \frac 12 p(p-1) \int_0^t  |u^n_{s\wedge \tau^n_R}|_H^{p-2}
   \|G(u^n_{s\wedge \tau^n_R})\|^2_{HS(H)} 1_{\{s < \tau^n_R\}} ds
\end{align*}
By assumption \eqref{crescG}
\begin{equation}\label{prima stima}
\begin{split}
 |u^n_{t\wedge \tau^n_R}|_H^p 
& +
    p\nu \int_0^t |u^n_{s\wedge \tau^n_R}|_H^{p-2} 
         \|u^n_{s\wedge \tau^n_R}\|_V^2 1_{\{s <\tau^n_R\}}ds
\\
&\le
 \left|  u_0^n\right|_H^p 
  + p   \left|\widetilde{M}_{t}^{n}\right|
  +\frac 12 p(p-1) \int_0^t  |u^n_{s\wedge \tau^n_R}|_H^{p-2}
   \left(\lambda_0|u^n_{s\wedge \tau^n_R}|^2_H+\rho\right) 
       1_{\{s <\tau^n_R\}}ds
\end{split}
\end{equation}
where
\[
 \widetilde{M}_{t}^{n}=
\int_0^t   |u^n_{s\wedge \tau^n_R}|_H^{p-2}
       \langle u^n_{s\wedge \tau^n_R},1_{\{s < \tau^n_R\}}
       \pi_n G(u^n_{s\wedge \tau^n_R}) dw(s) \rangle_H 
\]
is a square integrable martingale.
\\
Therefore
\begin{equation}\label{sup-noE}
\begin{split}
\sup_{t \in [0,r]} |  u^n_{t\wedge \tau^n_R}|_H^p
 \le 
&
\left|  u_0^n\right|_H^p 
  +p  \sup_{t \in [0,r]} \left|\widetilde{M}_{t}^{n}\right|
\\
& 
 +\frac 12 p(p-1) (\lambda_0+\rho) 
      \int_0^r  |u^n_{s\wedge \tau^n_R}|_H^p 1_{\{s<\tau^n_R\}}ds
  + \frac 12 p(p-1) \rho r
\end{split}
\end{equation}
By Burkholder-Davis-Gundy inequality, we estimate the supremum of the
martingale $ \widetilde{M}_{t}^{n}$;
for some constant $C>0$ we have
\begin{multline*}
 p E\sup_{0\le t\le r} \left| \int_0^t |u^n_{s\wedge \tau^n_R}|^{p-2}_H  
      \langle u^n_{s\wedge \tau^n_R},1_{\{s<\tau^n_R\}}
          \pi_nG(u^n_{s\wedge \tau^n_R})dw(s)\rangle_H \right| 
\\
\le
 C p E \left( \int_0^r |u^n_{s\wedge \tau^n_R}|^{2p-2}_H
              \|G(u^n_{s\wedge \tau^n_R})\|^2_{HS(H)} 1_{\{s<\tau^n_R\}}ds
                 \right)^{1/2}
\end{multline*}
Then by assumption (\ref{crescG}) 
\begin{equation}\label{BDG}
\begin{split}
p E \sup_{t \in [0,r]} & \left|\widetilde{M}_{t}^{n}\right|
\\
&\le E \left[ \sup_{0\le t\le r} |u^n_{t\wedge \tau^n_R}|_H^{p/2}
                 C p \left(\int_0^r  |u^n_{s\wedge \tau^n_R}|^{p-2}_H
     (\lambda_0 |u^n_{s\wedge \tau^n_R}|^2_H+\rho) 1_{\{s<\tau^n_R\}}ds \right)
          ^{1/2}\right]
\\
&\le
 \left(E \sup_{0\le t\le r} |u^n_{t\wedge \tau^n_R}|_H^p\right)^{1/2}
 \left(E C^2 p^2 \int_0^r  |u^n_{s\wedge \tau^n_R}|^{p-2}_H
                     (\lambda_0 |u^n_{s\wedge \tau^n_R}|^2_H+\rho)
                   1_{\{s<\tau^n_R\}}ds \right)^{1/2}
\\
&\le
 \frac 12 E  \sup_{0\le t\le r} |u^n_{t\wedge \tau^n_R}|_H^p
   + \frac 12 C^2 p^2  (\lambda_0 +\rho) E\int_0^r |u^n_{s\wedge \tau^n_R}|^p_H
      1_{\{s<\tau^n_R\}} ds +\frac 12 C^2 p^2 \rho r
\end{split}
\end{equation}
Then, by \eqref{sup-noE} and \eqref{BDG}
\begin{equation}
\begin{split}
\frac 12 E\sup_{t \in [0,r]} |  u^n_{t\wedge \tau^n_R}|_H^p
& \le 
  E|u_0^n|_H^p 
    + \frac 12 C^2 p^2  (\lambda_0 +\rho)
      \int_0^r  E|u^n_{s\wedge \tau^n_R}|_H^p 1_{\{s<\tau^n_R\}}ds
  + \frac 12 C^2 p^2 \rho r
\\
& \le
  E|u_0^n|_H^p 
    + \frac 12 C^2 p^2  (\lambda_0 +\rho)
      \int_0^r E\sup_{t \in [0,s]}|u^n_{t\wedge \tau^n_R}|_H^p ds
  + \frac 12 C^2 p^2 \rho r
\end{split}
\end{equation}
By Gronwall lemma, for any $r>0$ we have
\begin{equation}\label{mediap}
 E \sup_{0\le t\le r} |u^n_{t\wedge \tau^n_R}|_H^p \le C_1
\end{equation}
for some positive constant $C_1=C_1(p,T,\lambda_0,\rho,m_p)$ independent of
$n$ and $R$.
Here $m_p=E|u_0|_H^p$.
Notice that
$ E\left[  \left|  u_{0}^{n}\right|_{H}^p \right]  
 \leq m_p$.

Coming back to \eqref{prima stima}, with similar arguments
we also obtain
\[
 E \int_0^T |u^n_{s\wedge\tau_R}|_H^{p-2} 
              \|u^n_{s\wedge\tau_R}\|_V^2 1_{\{s< \tau_R\}}ds
\le C \qquad \forall n,R
\]
for a new positive constant $C$ depending on $m_p, p, \lambda_0,
\rho,T$ but not on $n,R$. For $p=2$ we have
\begin{equation}\label{stimaV}
 E \int_0^T \|u^n_{s\wedge\tau_R}\|_V^2 1_{\{s< \tau_R\}}ds
\le C_2
\end{equation}
for some positive constant $C_2=C_2(T,\lambda_0,\rho,m_2)$ 
independent of $n$ and $R$.

For the last estimate, we proceed as follows. 
First, from \eqref{Ipo-phi} we have that
\begin{align*}
 \int_{\mathcal T}|u(x)|^{1+b}dx
& =
  \int_{\mathcal T}|u(x)|^{1+b} 1_{\{|u(x)|\le K\}} dx+
  \int_{\mathcal T}|u(x)|^{1+b} 1_{\{|u(x)| > K\}} dx
\\
& \le
  K^{1+b}|\mathcal T| + \int_{\mathcal T} |u(x)|^2 \frac 1 {a_1}
  \sigma(|u(x)|) dx
\\
& =
 K^{1+b}|\mathcal T| + \frac 1 {a_1}\langle \Phi(u),u\rangle_H
\end{align*}
i.e.
\begin{equation} \label{stime5}
 \langle \Phi(u),u\rangle_H 
  \ge 
  a_1 |u|^{1+b}_{\mathbb L^{1+b}} -a_1 K^{1+b}
 |\mathcal T|
\end{equation}
Moreover
\begin{equation}\label{BA-1}
 2 |\langle B(u,u),A^{-1}u\rangle|
 \le
 2 |u|^2_{\mathbb L^4} |u|_{V^\prime}
 \le 
  |u|^4_{\mathbb L^4} + |u|^2_{V^\prime}
 \le
 a_1 |u|^{1+b}_{\mathbb L^{1+b}}+ C +C|u|^2_H
\end{equation}
Apply It\^{o} formula to $|u^n_{t\wedge \tau^n_R}|_{V^\prime}^2
=\langle u^n_{t\wedge \tau^n_R}, A^{-1} u^n_{t\wedge \tau^n_R} \rangle_H$
and get
\begin{align*}
  |u^n_{t\wedge \tau^n_R}|_{V^\prime}^2 
   =
&  |u_0^n|_{V^\prime}^2
  -2\int_0^t \langle \Phi (u^n_{s\wedge \tau^n_R}),
         u^n_{s\wedge \tau^n_R}\rangle_H 1_{\{s<\tau^n_R\}}ds
\\
 & - 2\int_0^t\langle \pi_n B(u^n_{s\wedge \tau^n_R},u^n_{s\wedge
   \tau^n_R}),
         A^{-1}u^n_{s\wedge \tau^n_R}\rangle_H 1_{\{s<\tau^n_R\}}ds
\\
 &  +\underline{M}_{t}^{n}+\int_0^t \|\pi_n
   A^{-1/2}G(u^n_{s\wedge \tau^n_R})\|_{HS(H)}^2 1_{\{s<\tau^n_R\}}ds
\end{align*}
where 
\[
 \underline{M}_{t}^{n}
 =
 2\int_0^t \langle u^n_{s\wedge \tau^n_R},1_{\{s<\tau^n_R\}}
      \pi_n A^{-1}G(u^n_{s\wedge \tau^n_R}) dw(s)\rangle_H
\]
is a square integrable martingale.
\\
Use \eqref{stime5}  and \eqref{BA-1}; then
\begin{align*}
 E |u^n_{T\wedge \tau^n_R}|_{V^\prime}^2 
 &+
  a_1E \int_0^T |u^n_{s\wedge \tau^n_R}|^{1+b}_{\mathbb L^{1+b}} 
          1_{\{s<\tau^n_R\}} ds
\\
& \le
  E |u_0^n|_{V^\prime}^2
 + CT E \sup_{0\le t\le T} |u^n_{t\wedge \tau^n_R}|_H^2 
\\
& \qquad 
 + C \lambda_0 E \int_0^T|u^n_{s\wedge \tau^n_R}|^2_H
 1_{\{s<\tau^n_R\}} ds+C (1+\rho) T 
\end{align*}
According to \eqref{mediap} (for $p=2$), we conclude that
\begin{equation}\label{stima5}
 E \int_0^T |u^n_{s\wedge \tau^n_R}|^{1+b}_{\mathbb L^{1+b}}
 1_{\{s<\tau^n_R\}} ds
 \le C_3
\end{equation}
for some positive constant $C_3=C_3(T,\lambda_0,a_1,C_1,C_2)$. 
          \vspace{1mm}

\section{Appendix 2: a compactness result}
Let $\mathbb X$ be the closure of $\mathcal D^\infty$ 
w.r.t. the norm 
$$
 |u|_{\mathbb X}:=
  \left(\int_{\mathcal T}\{|u(x)|^4 |\nabla u(x)|^2 +4 |u(x)|^2
  \sum_{i=1}^3 [u(x)\cdot \partial _i u(x)]^2 \}dx
  \right)^{1/6}
$$
\begin{theorem} \label{compatt}
$\mathbb X \subset\mathbb{L}^{6}\cap H$ and the immersion is compact.
\end{theorem}

\begin{proof}
First observe that for smooth fields $u$ we have
\[
 \partial_i \left(|u|^2 u\right)  
 =
 |u|^2 \partial_i u +2u (u\cdot \partial_i u)
\]
Hence
\[
 \left\||u|^2 u\right\|_V^2
 \le 
 C \int_{\mathcal{T} }|u|^{4}|\nabla u|^2 dx
 \le 
 C |u|^6_{\mathbb X}
\]
and thus by Poincar\'{e} inequality
\[
 |u|_{\mathbb L^6}^{6}
 =
 \left| |u|^2 u\right|_{\mathbb L^2}^2
 \le 
  \left\| |u|^2 u\right\|^2_V
 \le
 C^{\prime} |u|^6_{\mathbb X}
\]
This proves that the closure of $\mathcal D^\infty$ with respect to
the $\mathbb L^6$-norm is a space bigger than its closure with respect to the
$\mathbb X$-norm;
hence
$\mathbb X \subset\mathbb{L}^{6}\cap H$.

Moreover, if $\{u_n\}$ is a bounded sequence in
$\mathbb X$, it is
bounded in $\mathbb{L}^{6}$ and $\{\left\|  \left|  u_{n}\right|  ^{2}
u_{n}\right\|  _{V}^{2}\}$ is also bounded. By Rellich Theorem, the sequence
$\{  \left|  u_{n}\right|  ^{2}u_{n}\}$ is relatively compact in
$\mathbb{L}^{2}$ and so
there exists  a subsequence $\{  \left|  u_{n_{k}}\right|^2 u_{n_k}\}  $
converging strongly in $\mathbb{L}^{2}$ to some field $\xi$; we also
have that 
$\{u_{n_k}\}$ converges weakly in $\mathbb{L}^{6}$
to some field $u$. The strong convergence implies in particular that
\[
\left|  u_{n_{k}}\right|  _{\mathbb L^{6}}^{6}
 =\left|  \left|  u_{n_{k}}\right|^{2}u_{n_{k}}\right|  _{\mathbb
 L^{2}}^{2}\rightarrow\left|  \xi\right|  _{\mathbb L^{2}}^{2}.
\]
Thus, if we prove that 
$\left|  \xi\right|  _{\mathbb L^{2}}^{2}=\left|  u\right|
_{\mathbb L^{6}}^{6}$, then from the weak convergence of $\{u_{n_k}\}$
to $u$ in $\mathbb{L}^{6}$ and the convergence of norms $\left|  u_{n_{k}
}\right|  _{\mathbb L^{6}}^{6}\rightarrow\left|  u\right|  _{\mathbb L^{6}}^{6}$, we deduce
that $\{u_{n_k}\}$ converges strongly to $u$ in $\mathbb{L}
^{6}$ and the proof of the compact embedding will be complete.

So it remains to show that $\left|  \xi\right|  _{\mathbb L^{2}}^{2}=\left|  u\right|
_{\mathbb L^{6}}^{6}$. Let us introduce the function
\[
v\left(  x\right)  =\left\{
\begin{array}
[c]{ccc}%
\frac{\xi\left(  x\right)  }{\left|  \xi\left(  x\right)  \right|  ^{2/3}} &
\text{if} & \xi\left(  x\right)  \neq0\\
0 & \text{if} & \xi\left(  x\right)  =0
\end{array}
\right.  .
\]
Let us prove that there is a subsequence $\{ u_{n_k^\prime}\}$
such that
\begin{equation}
 u_{n_k^{\prime}}  \rightarrow v\text{ a.s. on }%
\mathcal{T}.\label{a.s. convergence}%
\end{equation}
This implies $v=u$ (the a.s. limit and the $\mathbb{L}^{6}$ weak limit must
coincide, since by Vitali theorem there is strong convergence in any
$\mathbb{L}^{p}$ with $p<6$). Since $\left|  v\left(  x\right)  \right|
^{6}=\left|  \xi\left(  x\right)  \right|  ^{2}$ where $\xi\left(  x\right)
\neq0$, we have $\left|  \xi\right|  _{\mathbb L^{2}}^{2}=\left|  v\right|  _{\mathbb L^{6}%
}^{6}=\left|  u\right|  _{\mathbb L^{6}}^{6}$, as we want.

Thus it remains to prove (\ref{a.s. convergence}). The strong convergence
above implies that there is a subsequence 
$\{ |u_{n_k^\prime}|^2 u_{n_k^\prime}\}  $ that converges to $\xi$ a.s. on
$\mathcal{T}$. Let $x\in\mathcal{T}$ be such that 
$|u_{n_k^\prime}(x)|^2 u_{n_k^\prime}(x)   \rightarrow
\xi(x)$. Taking the norm in
$\mathbb{R}^{3}$, this implies that $|u_{n_k^\prime}(x)|^3 \rightarrow
|\xi(x)|$, 
hence $|u_{n_k^\prime}(x)|  \rightarrow
|\xi(x)|^{1/3}$. If $\xi(x)  =0$,
this implies $u_{n_k^\prime}(x)  \rightarrow 0$, as we want in
(\ref{a.s. convergence}). If $\xi\left(  x\right)  \neq0$, this implies
$\left|  u_{n_{k}^{\prime}}\left(  x\right)  \right|  \neq0$ eventually and
\[
u_{n_{k}^{\prime}}\left(  x\right)  =\frac{\left|  u_{n_{k}^{\prime}}\left(
x\right)  \right|  ^{2}u_{n_{k}^{\prime}}\left(  x\right)  }{\left|
u_{n_{k}^{\prime}}\left(  x\right)  \right|  ^{2}}\rightarrow\frac{\xi\left(
x\right)  }{\left|  \xi\left(  x\right)  \right|  ^{2/3}}.
\]
Thus (\ref{a.s. convergence}) is true. The proof is complete.
\end{proof}

\end{document}